\documentclass[1p,times]{elsarticle}
\usepackage{hyperref}
\usepackage{amssymb,amsthm}
\usepackage[textsize=small]{todonotes}
\usepackage{graphicx}
\usepackage{upref,setspace}
\usepackage{xcolor,colortbl}
\usepackage{enumerate}
\usepackage{bbm}
\usepackage{float}
\usepackage{bm}
\usepackage[mathscr]{eucal}
\usepackage{graphicx}
\usepackage{sidecap}
\usepackage{latexsym}
\usepackage{psfrag}
\usepackage{epsfig}
\usepackage{enumerate}
\usepackage{amsmath}
\usepackage{amsthm}
\usepackage{amssymb}
\usepackage{multirow}
\usepackage{comment}
\usepackage{amsfonts}
\usepackage{color}
\usepackage[textsize=small]{todonotes}
\usepackage{hyperref}
\usepackage{mathtools}

\usepackage[T1]{fontenc}

\newtheorem{theorem}{Theorem}

\newtheorem{lemma}[theorem]{Lemma}
\newtheorem{construction}[theorem]{Construction}

\newtheorem{assumption}[theorem]{Assumption}

\newtheorem{example}[theorem]{Example}
\newtheorem{remark}[theorem]{Remark}

\newcommand{\veps}{\varepsilon}

\newcommand{\PP}{\mathbb{P}}

\newcommand{\RR}{\mathbb{R}}

\newcommand{\NN}{\mathbb{N}}

\newcommand{\clv}{\mathcal{V}}
\newcommand{\clu}{\mathcal{U}}
\newcommand{\clm}{\mathcal{M}}

\newcommand{\clr}{\mathcal{R}}
\newcommand{\cla}{\mathcal{A}}
\newcommand{\clb}{\mathcal{B}}

\newcommand{\clf}{\mathcal{F}}

\newcommand{\clp}{\mathcal{P}}

\newcommand{\bfL}{\mathbf{L}}

\begin{document}

\begin{frontmatter}

\title{Empirical Measure Large Deviations for  Reinforced Chains on Finite Spaces}
\author[1]{Amarjit Budhiraja}
\ead{budhiraj@email.unc.edu}

\author[2]{Adam Waterbury}
\ead{awaterbury@ucsb.edu}
\address[1]{Department of Statistics and Operations Research, University of North Carolina at Chapel Hill,
Hanes Hall, Chapel Hill, NC 27599, USA}
\address[2]{Department of Statistics and Applied Probability, South Hall, University of California, Santa Barbara,
CA 93106, USA}

\begin{abstract}
Let $A$ be a transition probability kernel on a finite state space $\Delta^o =\{1, \ldots , d\}$ such that $A(x,y)>0$ for all $x,y \in \Delta^o$. Consider a reinforced chain given as a sequence $\{X_n, \; n \in \NN_0\}$ of $\Delta^o$-valued random variables,  defined recursively according to,
$$L^n = \frac{1}{n}\sum_{i=0}^{n-1} \delta_{X_i}, \;\; P(X_{n+1} \in \cdot \mid X_0, \ldots, X_n) = L^n A(\cdot).$$
We establish a large deviation principle for $\{L^n\}$. The rate function takes a strikingly different form than the Donsker-Varadhan rate function associated with the empirical measure of the Markov chain with transition kernel $A$ and is described in terms of a novel deterministic infinite horizon discounted cost control problem with an associated linear controlled dynamics and a nonlinear running cost involving the relative entropy function. Proofs are based on an analysis of time-reversal of controlled dynamics in  representations for log-transforms of exponential moments, and on weak convergence methods.
\end{abstract}


\end{frontmatter}

\section{Introduction}

Processes with reinforced dynamics have been used to model systems in ecology and biology (see, e.g., \cite{Donnelly1986PartitionSP, hop1}) and have been applied in a wide range of sampling and optimization problems, including non-linear Markov chain Monte Carlo \cite{andjasdou} and stochastic optimization \cite{andtadvlavih, biaforhac, for1}. In many such settings, the time-asymptotic properties of processes have been extensively studied, and laws of large numbers and central limit theorems have been established under broad conditions. In the current work our focus is on the study of large deviation asymptotics for certain types of reinforced dynamics on a finite state space which, to the best of our knowledge, have not been studied in the literature to date. One of the reasons for a lack of results in this direction is that the state dynamics in such systems is not Markovian and the evolution of the state depends on the full path history, making the usual methods of empirical measure large deviation analysis challenging to implement. To motivate the problem of interest consider first the elementary setting of a Markov chain $\{X^0_n, \; n \in \NN_0\}$ with values in the finite state space $\Delta^o = \{1, \ldots , d\}$ and a transition probability matrix $A$ with $A(x,y)>0$ for all $x,y \in \Delta^o$. The celebrated results of Donsker and Varadhan \cite{donvar1, donvar3} give a large deviation principle (LDP) for the empirical measure sequence 
$\{L^n_0, \; n \in \NN\}$ defined as $L^n_0 = n^{-1} \sum\limits_{i=0}^{n-1} \delta_{X^0_i}$, $n \in \NN$, with rate function
given as 
\begin{equation}\label{eq:830n}
I(\theta) = \inf_{\gamma \in \cla(\theta)} R(\gamma \| \theta \otimes A), \; \theta \in \clp(\Delta^o),
\end{equation}
where $\clp(\Delta^o)$ is the space of probability measures on $\Delta^o$, $\theta \otimes A$ is a probability measure on
$\Delta^o\times \Delta^o$ defined as, for $x,y \in \Delta^o$, $\theta \otimes A(x,y) \doteq  \theta(x) A(x,y)$,  $\cla(\theta)$ is the space of all probability measures $\gamma$ on $\Delta^o\times \Delta^o$ for which the two marginals are the same as $\theta$, i.e.,
$\sum_{y \in \Delta^o} \gamma(x,y) = \sum_{y \in \Delta^o} \gamma(y,x) = \theta(x)$, for all $x \in \Delta^o$, and $R(\cdot \|\cdot)$ is the relative entropy function (see Section \ref{sec:not}).
Consider now the reinforced chain $\{X_n, \; n \in \NN_0\}$ associated with the  transition probability kernel $A$ which is constructed recursively as $X_0 = x_0$ for some $x_0 \in \Delta^o$ and, having defined $X_0, \ldots , X_{n-1}$
and $L^n = n^{-1}\sum\limits_{i=0}^{n-1} \delta_{X_i}$, the conditional distribution of $X_n$ given $\{X_0, \ldots , X_{n-1}\}$
is $L^nA(\cdot) = \sum_{x\in \Delta^o} L^n(x)A(x, \cdot)$. Thus, at each time instant, one of the previously visited sites $x^*$ is chosen at random (with probabilities proportional to visit frequencies) and then the new site is selected according to the distribution $A(x^*, \cdot)$. The law of the large numbers for $L^n$ and $L^n_0$ is the same, namely both converge to the unique stationary distribution of the Markov chain $\{X^0_n\}$. However, as we will see, the study of the large deviation behavior of $\{L^n, \; n \in \NN\}$ requires a rather different type of analysis than the sequence $\{L^n_0, \; n \in \NN\}$, and the associated rate function, which is introduced in \eqref{def:ratefunction1}, has a strikingly different form than the rate function in \eqref{eq:830n}. In particular, the emergence of a discount factor in the rate function and the rate function's characterization as the value function of an infinite horizon discounted cost problem is novel in the context of empirical measure large deviation theory.

One special case of our results is the following. Let $p = (p_x)_{x \in \Delta^o}$ be a  positive probability vector on $\Delta^o$. Then Sanov's theorem (cf. \cite[Theorem 3.3]{buddupbook}) tells us that if $\{X^0_n\}$ is an iid sequence with law $p$, then the empirical measure $L^n_0 = n^{-1} \sum\limits_{i=0}^{n-1} \delta_{X^0_i}$ satisfies a LDP with rate function 
$I(\theta) = R(\theta \| p)$, $\theta \in \clp(\Delta^o)$. Consider now the setting where $p^0$ is a positive probability vector on $\Delta^o \cup \{0\}$ and $p_x \doteq p^0_x/(1-p^0_0)$, $x \in \Delta^o$. We construct a chain $\{X_n\}$ for which the conditional distribution of $X_n$ given $\{X_0, \ldots, X_{n-1}\}$ is $p^n$, where $p^n_x = p^0_x + (1- p^0_0)L^n$, $x \in \Delta^o$, where
$L^n = n^{-1}\sum\limits_{i=0}^{n-1} \delta_{X_i}$. For this chain, at each instant the new state is chosen, independently of the past, according to the probability vector $p^0$, but if the chosen state is $0$, the chain immediately moves to a state chosen at random from the collection of previously visited states (with probabilities proportional to visit frequencies).
These types of reinforcement mechanisms have been used for numerical approximations of quasi-stationary distributions of Markov chains (see, e.g., \cite{AFP88, benclo}).
Once again the law of large numbers for $L^n_0$ and $L^n$ are the same, namely both converge to the probability measure $p = p^0/(1-p^0_0)$, however, as will be seen, the large deviation behavior of $L^n$ is more complex and is governed by an infinite horizon discounted cost problem. Another special case is where $A(x,y) = p_y$, $x,y \in \Delta ^o$, where $p$ is a probability vector on $\Delta^o$. In this case the model simply reduces to an iid sequence (with no reinforcement) for which the empirical measure large deviation principle is given by Sanov's theorem. In this special case the rate function given in \eqref{def:ratefunction1} is easily seen to be the same as the familiar rate function in Sanov's theorem (see Remark \ref{rem:sanov}).

We now comment on proof techniques. The starting point is a reformulation of the large deviation principle in terms of Laplace asymptotics (cf. \cite[Theorem 1.8]{buddupbook}) and a stochastic control representation using relaxed controls (cf. \eqref{eq:varrep}) for log-transforms of exponential moments of empirical measure functionals. 
Inspired by the ODE-method for the study of stochastic approximation schemes (cf. \cite{BMP12, bor09, kus03}), we introduce a suitable  continuous time interpolation for the controlled processes in the stochastic control representation.   The impact of reinforcement on the large deviation behavior becomes evident when one considers the asymptotics of these interpolated controlled sequences. It turns out that this asymptotic behavior is particularly well understood when one considers the time-reversed trajectories (for interpolated controlled sequences) going back from states far off in the future. The asymptotic time-reversed paths have a simple form linear dynamics which also reveals the exponential discounting of the contribution, to the cost, of initial segments of near-optimal paths. Proof of the large deviation upper bound relies on tightness and characterization of weak limits of time-reversed trajectories and their costs. For the lower bound,  one takes a constructive approach. We first argue that by suitable approximations, mollification, and discretization one can find simple form piecewise-constant near-optimal trajectories for the variational problem describing the Laplace asymptotics. The remaining work is to then construct suitable controlled state sequences for which the associated interpolated paths and the corresponding costs converge in probability to these simple form near-optimal paths and associated costs.

We remark that the current work considers the simplest forms of reinforcement mechanisms in finite-state models. A general form of reinforced dynamics corresponds to a setting where the conditional distribution of $X_{n+1}$ given $\{X_0, \ldots, X_{n-1}\}$ is given as $L^n A(L^n)$, where $A$ is a suitable (nonlinear) map from $\clp(\Delta^o)$ to the space of stochastic kernels on $\Delta^o \times \Delta^o$ and $L^n$ is as before the empirical measure. This general setting introduces significant new challenges, particularly in the treatment of the lower bound, and is a topic currently under study.

\subsection{Notation.} \label{sec:not} The following notation is used.
Fix $d \in \NN$.
Let $\Delta \doteq \{0, 1, \ldots, d\}$, and let $\Delta^o = \{1, \ldots , d\}$. For a metric space $S$, $\clb(S)$ denotes the corresponding Borel $\sigma$-field and $\clp(S)$ denotes the space of probability measures on $(S, \clb(S))$ equipped with the topology of weak convergence. Recall that a function $I: S \to [0, \infty]$ is called a rate function if it has compact sublevel sets, namely
$S_k \doteq \{ x \in S: I(x) \le k\}$ is compact for every $k \in [0, \infty)$. For $x \in S$, $\delta_x \in \clp(S)$ is the Dirac probability measure concentrated at the point $x$. 
For $\nu, \mu \in \clp(S)$, we denote the relative entropy of $\nu$ with respect to $\mu$
as $R(\nu \|\mu)$, which is the extended real number defined as
\[
R(\nu \|\mu) \doteq \int_S \left(\log \frac{d\nu}{d\mu}\right) d\nu,
\]
if $\nu$ is absolutely continuous with respect to $\mu$, and $+\infty$ otherwise.
Let $\clv^d \doteq \{e_1, \ldots , e_d\}$, where $e_x$ is the $x$-th unit coordinate vector in $\RR^d$, and denote by $\clm(\clv^d\times \RR_+)$ the space of locally finite measures on $\clv^d\times \RR_+$ with the vague topology. 
Denote by $\clu$ the collection of all measurable maps from $\RR_+$ to $\clp(\Delta^o)$. We denote by $C_b(\clp(\Delta^o))$  the space of bounded continuous functions from $\clp(\Delta^o)$ to $\RR$. For $m, \tilde m \in \clp (\Delta^o)$, we write  $\| m - \tilde m\| \doteq \sum\limits_{x\in\Delta^o} |m(x) - \tilde m(x)|$. We use the same notation for the norm of a vector in $\RR^d$. As a convention $\int_a^b f(s) ds$ is taken to be $0$ if $a\ge b$. For $\RR^d$-valued random variables $\{X_n, n \in \NN\}, X$, we say that $X_n \to X$ in $\mathcal{L}^2$ as $n \to \infty$ if $E\| X_n - X\|^2 \to 0$ as $n\to\infty$. For a vector $v \in \RR^d$ we use the notation $v_x$ and $v(x)$ interchangeably to denote the $x$-th coordinate of $v$.

\subsection{Description of the Model}

Consider a map $K : \clp(\Delta^o) \to \clp(\Delta^o)$ satisfying Assumption \ref{assu:kassumptions} below.

\begin{assumption}\label{assu:kassumptions}
There is a $d \times d$ stochastic matrix $A$ such that $\delta_0 \doteq \inf_{x,y \in \Delta^o} A_{x,y} > 0$, and, for $m \in \clp(\Delta^o)$, $K(m) = m A$.
\end{assumption}

For fixed $x_0 \in \Delta^o$, we consider a collection $\{X_n, \; n \in \NN_0\}$ of $\Delta^o$-valued random variables, a collection $\{L^n,\;  n \in \NN\}$ of $\mathcal{P}(\Delta^o)$-valued random measures, and a filtration $\{\mathcal{F}_n, \; n \in \NN_0\}$ on some probability space $(\Omega,\mathcal{F}, \PP)$, defined recursively as follows. Let $X_0  \doteq x_0$, $\mathcal{F}_0 \doteq \{\emptyset, \Omega\}$,  and $L^1 \doteq \delta_{x_0}$.  Having defined 
$\{X_i, L^{i+1},\;  0\le i \le n\}$
and $\sigma$-fields $\{\clf_i, \; i \le n\}$ for some $n \in \NN_0$, define $\PP(X_{n+1} = y | \mathcal{F}_n) \doteq K(L^{n+1})(y)$, $y \in \Delta^o$,
 $\mathcal{F}_{n+1} \doteq \sigma\{X_k, \; k \leq n+1\}$,
and
\begin{equation}\label{eq:occmzr}
L^{n+2} \doteq \frac{1}{n+2} \sum\limits_{i=0}^{n+1} \delta_{X_i}.
\end{equation}

For $m \in \clp( \Delta^o)$, let $\rho(m) \in \clp(\clv^d)$ be defined as
\[
\rho(m)(e_x) \doteq K(m)(x) = (mA)_x, \quad x \in \Delta^o.
\]
Let $\{\nu^k(m), \; m\in \clp(\Delta^o)\}_{k\in \NN}$ be iid random fields with values in $\clv^d$ such that, for each $m \in \clp(\Delta^o)$,
\begin{align*}
	 \PP(\nu^1(m) = e_x) = \rho(m)(e_x) = K(m)(x), \quad x \in \Delta^o.
\end{align*}
Then, we can write the evolution equation for $L^n$ as 
\begin{equation}
	L^{n+1} = L^n + \frac{1}{n+1}\left[\nu^n(L^n) - L^n\right],\quad n\in \NN.
\end{equation}

\subsection{Statement of Results}

For $m \in \clp(\Delta^o)$, let  $\clu(m)$ be the collection of all $\eta \in \clu$ such that if $M : \RR_+ \to \clp(\Delta^o)$ satisfies 
\begin{equation}\label{eq:Meq}
	M(t) = m - \int_0^t \eta(s) ds + \int_0^t M(s) ds, \;\; t \in \RR_+, 
	\end{equation}
then $M \in C(\RR_+: \clp(\Delta^o))$. 
Note that, given $\eta \in \clu$ and $m \in \clp(\Delta^o)$, \eqref{eq:Meq} always has a unique solution in $C(\RR_+: \RR^d)$; however for such an $\eta$ to be in $\clu(m)$ we require that the solution is in fact in $C(\RR_+: \clp(\Delta^o))$.

Given $M \in C(\RR_+: \clp(\Delta^o))$ that satisfies  \eqref{eq:Meq}, with some
$m \in \clp(\Delta^o)$ and $\eta \in \clu$, for all $t \in \RR_+$, we say that $M$ solves $\clu(m,\eta)$.

Define $I: \clp(\Delta^o) \to \RR_+$ as
\begin{equation}\label{def:ratefunction1}
	I(m) \doteq \inf_{\eta \in \clu(m)} \int_0^{\infty} \exp(-s) R\left(\eta(\cdot \mid s) \| K(M(s))\right) ds,\;\;  m \in \clp(\Delta^o),
\end{equation}
where $M$ solves $\clu(m,\eta)$.
For a $\eta \in \clu$, define $\Lambda^{\eta}: \RR_+ \to \clp(\clv^d)$ as
$\Lambda^{\eta}(s)(e_x) \doteq \eta(s)(x)$, for $s\in \RR_+$ and $x \in \Delta^o$.
Then the above rate function can equivalently be written as
\begin{equation}\label{eq:equivrf}
	I(m) \doteq \inf_{\eta \in \clu(m)} \int_0^{\infty} \exp(-s) R\left(\Lambda^{\eta}(\cdot \mid s) \| \rho(M(s))\right) ds,\;\;  m \in \clp(\Delta^o),
\end{equation}
where $M$ solves $\clu(m,\eta)$.

The following theorem is the main result of this work, which establishes a large deviations principle (LDP) for  $\{L^{n}, \;  n \in \NN\}$.
\begin{theorem}\label{thm:ldp1}
Let $I : \clp(\Delta^o) \to  \RR_+$ be the  function defined in \eqref{def:ratefunction1}. Then $I$ is a rate function and the sequence $\{L^{n+1},\; n \in \NN\}$ satisfies an LDP with rate function $I$. Namely, for each closed set $F \subseteq \clp (\Delta^o)$,
\[
\limsup_{n\to\infty}n^{-1} P(L^{n+1} \in F) \le -\inf_{m \in F}I(m),
\]
and for each open set $G \subseteq \clp(\Delta^o)$,
\[
\liminf_{n\to\infty} n^{-1}P(L^{n+1} \in G) \ge - \inf_{m \in G}I(m).
\]
\end{theorem}
\begin{proof}
In view of \cite[Theorem 1.8]{buddupbook} it suffices to show that
for every $F \in C_b(\clp(\Delta^o))$,
	\begin{equation*}
	\liminf_{n\to \infty}-n^{-1} \log E \exp[-n F(L^{n+1})] = \inf_{m \in \clp(\Delta^o)} [F(m) + I(m)]
	\end{equation*}
	and that $I$ is a rate function (namely it has compact sublevel sets).
	The first  statement is shown in Theorems \ref{thm:laplaceupperbound1} and \ref{thm:lowbd}, which establish the Laplace upper bound and Laplace lower bound respectively, whereas the second statement is shown in Section \ref{sec:levelset}.
\end{proof}

\subsection{Examples}

We note three examples that are covered by the model studied in this work.

\begin{example}
\begin{enumerate}[(a)]
\item Let $p \in \mathcal{P}(\Delta)$ satisfy $\inf_{x\in\Delta} p_x > 0$, and let $p^o$ be the sub-probability  measure obtained by restricting $p$ to $\Delta^o$. Consider the transition matrix $A$ given by  $A = P^o + p_0 I,$ where $P^o$ is the $d\times d$ matrix whose every row is $p^o$, and $I$ is the $d\times d$ identity matrix. Then, Assumption \ref{assu:kassumptions} is satisfied and the corresponding $\{L^{n+1}, \; n \in \NN_0\}$ is a special case of the algorithm for approximating quasi-stationary distributions studied in \cite{AFP88}, for the case when the system's dynamics are not state-dependent. 

\item Let $A$ be a transition probability matrix on $\Delta^o$ satisfying $\inf\nolimits_{x,y \in \Delta ^o}A_{x,y} > 0$. Then, Assumption \ref{assu:kassumptions} is satisfied and $\{X_n\}$ can be interpreted as a genetic-type algorithm for a population with uniform fitness, as $P(X_n = x | \clf_{n-1}) = \sum\limits_{i=0}^{n-1} n^{-1} A_{X_i,x}$, see \citep{delmormic1}. 

\item Fix $\alpha \in (0,1)$, let $p^o\in \mathcal{P}(\Delta^o)$ satisfy $\inf\nolimits_{x\in\Delta^o}p^o_x > 0$. Let $B$ be a transition probability matrix on $\Delta^o$. Then $A = \alpha p^o + (1-\alpha) B$ satisfies Assumption \ref{assu:kassumptions}
and describes a setting where, at each step, a new state is chosen according to the iid law $p^o$ with probability $\alpha$, and according to $L^n B$ with probability $1-\alpha$.
As $\alpha$ ranges from $0$ to $1$, these models interpolate between the reinforced setting of part (b) and the iid setting of Sanov's theorem (see also Remark \ref{rem:sanov}).
 \end{enumerate}
\end{example}
\begin{remark}\label{rem:sanov}
The setting of Sanov's theorem corresponds to the case  $A= P^o$ 
where, for some $p^o \in \clp(\Delta^o)$ satisfying $\inf\nolimits_{x\in \Delta^o}p_x > 0$, $P^o$ is the $d\times d$ matrix whose every row is $p^o$. In this case it is easy to verify that for $m \in \clp(\Delta^o)$ the infimum on the right side of \eqref{def:ratefunction1} is achieved at the constant function $\eta(s) = m$ for all $s\in \RR_+$ and, furthermore, with this choice of $\eta$, $M(s) = m$ and $\rho(M(s)) = p^o$ for all $s\in\RR_+$. Thus, the rate function in \eqref{def:ratefunction1} reduces to the familiar rate function
$R(m \|p)$ in Sanov's theorem.
\end{remark}

\subsection{Organization}
Rest of the paper is organized as follows. 
In Section \ref{sec:controlproc} we present the  stochastic control representation that is key in the proofs of both upper and lower bounds.
This section also presents some basic tightness and limit point characterization results. Section \ref{sec:laplaceupperbound} is devoted to the proof of the Laplace upper bound (Theorem \ref{thm:laplaceupperbound1}) while Section \ref{sec:lowbd} proves the 
Laplace lower bound (Theorem \ref{thm:lowbd}). Finally, Section \ref{sec:levelset} shows that $I$ has compact sublevel sets.

\section{A Stochastic Control Representation}\label{sec:controlproc}
The key ingredient in the proof is a certain stochastic control representation for exponential moments of functionals of the empirical measures $\{L^n,\; n \in \NN\}$, which we now present.

The controlled stochastic system, for each $n \in \NN$, is a sequence $\{ \bar{L}^{n,k} , k \in \NN\}$ of $\clp(\Delta^o)$-valued random variables which is defined recursively in terms of a collection of random probability measures on $\clv^d$, $\{\bar \mu^{n,k}, k \in \NN\}$,
where for each $k \in \NN$, $\bar \mu^{n,k}$ is $ \bar \clf^{n,k} \doteq \sigma(\{\bar L^{n,j}, 1\le j \le k\})$ measurable, and, having defined $\{ \bar{L}^{n,j}, 1\le j \le k\}$,
$\bar L^{n,k+1}$ is defined as
\begin{equation}\label{eq:barlnk}
	\bar L^{n,k+1} = \begin{cases}
	\delta_{x_0}, & k =0\\
	\bar L^{n,k} +\frac{1}{k+1}\left[\bar \nu^{n,k}- \bar L^{n,k}\right],& k \geq 1,\\
	\end{cases}
\end{equation}
where $\bar \nu^{n,k}$ is a $\clv^d$-valued random variable such that
$$P[\bar \nu^{n,k} = e_x \mid \bar \clf^{n,k}] = \bar \mu^{n,k}(e_x), \; x \in \Delta^o.$$
We denote the collection of all such {\em control} sequences $\{\bar \mu^{n,k}, k \in \NN\}$
as $\Theta^n$.
It is convenient to consider a continuous time interpolation of the sequences 
$\{\bar L^{n,k}, \; k \in \NN\}$. Define the time interpolation sequence $\{t_k, \; k \in \NN_0\}$ by $t_0 \doteq 0$, and $t_k = \sum_{j=1}^k (j+1)^{-1}$, for $k \geq 1$.
For each $n \in \NN$, define the $C(\RR_+ : \clp(\Delta^o))$-valued random variable $\bar{L}^n$ by linear interpolation: namely, for each $k \in \NN_0$,
\begin{equation}\label{eq:411}\bar L^n(t) \doteq 
\begin{cases}
\bar L^{n,k+1}, &\quad  t = t_k\\
\bar L^{n, k+1} +(k+2) (t- t_k)  [\bar L^{n, k+2} - \bar L^{n,k+1}], &t\in (t_k,t_{k+1}).
\end{cases}
\end{equation}
 
Consider random measures on $\clv^d \times [0,t_n]$ defined as follows: for $A \subseteq \clv^d$ and $B \in \clb[0,t_n]$,
\begin{equation}\label{eq:519}
\begin{aligned}
\bar\Lambda^n(A\times B) \doteq \int_B \bar\Lambda^n(A\mid t) dt, \;\; &
	\bar\mu^n(A\times B) \doteq \int_B \bar \mu^n(A \mid t) dt,
\end{aligned}
\end{equation}
where, for $k \le n-1$ and $t \in [t_k, t_{k+1})$,
\begin{align*}
\bar \Lambda^n(\cdot \mid t) \doteq \delta_{\bar \nu^{n,k+1}}(\cdot),\quad & 	\bar \mu^n(\cdot \mid t) \doteq \bar \mu^{n,k+1}(\cdot).
\end{align*}
The following variational representation follows from \cite[Theorem 4.5]{buddupbook}. For each $F \in  C_b(\clp(\Delta^o))$,
\begin{equation}\label{eq:varrep}
	-n^{-1}\log E \exp[-n F(L^{n+1})]
	= \inf_{\{\bar \mu^{n,i}\} \in \Theta^n}E \left[F(\bar L^n(t_n) ) + n^{-1} \sum\limits_{k=0}^{n-1} R\left(\bar \mu^{n,k+1} \| \rho(\bar L^{n,k+1})\right)\right].
\end{equation}
We now rewrite the right side above using the continuous time interpolation introduced in
\eqref{eq:411} and the `relaxed control' representation in \eqref{eq:519}.
For each $s \in \RR_+$, let $m(s) \doteq \sup\{ k : t_k \leq s\}$
and $a(s) \doteq t_{m(s)}$. It follows from the definitions above that, for $0\le t \le t_n$,
\begin{equation}\label{eq:522}
\begin{split}
	\bar L^n(t) &= \bar L^n(0) + \int_0^t \sum\limits_{v 
	\in \clv^d} (v-\bar L^n(a(s))) \bar \Lambda^n(v\mid s) ds.
	\end{split}
\end{equation}
Define $\psi_e: \RR_+ \to  \{2,3,\dots\}$ as $\psi_e(t_k) \doteq k+2$ for $k \in \NN_0$, and, for  $k \in \NN_0$ and $t \in [t_k, t_{k+1})$, define $\psi_e(t)$ by constant interpolation
as $\psi_e(t) \doteq \psi_e(t_k)$. From this definition it follows that
\begin{equation}
	n^{-1} \sum\limits_{k=0}^{n-1} R\left(\bar \mu^{n,k+1} \| \rho(\bar L^{n,k+1})\right) 
	=   n^{-1} \int_0^{t_n} \psi_e(s) R\left(\bar \mu^n(\cdot \mid s) \|  \rho(\bar L^n(a(s)))\right) ds. \label{eq:540}
\end{equation}
Define $\clp(\clv^d \times \RR_+)$-valued random variables as follows: for $t\in \RR_+$ and $x \in \Delta^o$, let
\begin{equation}\label{eq:525}
\begin{aligned}
\gamma^n(\{e_x\} \times [0,t]) &\doteq 	n^{-1} \int_{0}^{t_n \wedge t} \psi_e(t_n-s) \bar \Lambda^n(e_x \mid t_n-s) ds\\
\beta^n(\{e_x\} \times [0,t])&\doteq n^{-1} \int_{0}^{t_n \wedge t} \psi_e(t_n-s)\bar \mu^n(e_x  \mid t_n-s) ds\\
\theta^n(\{e_x\} \times [0,t]) &\doteq n^{-1} \int_{0}^{t_n \wedge t}  \psi_e(t_n-s) \rho(\bar L^n(a(t_n-s)))(e_x ) ds.
\end{aligned}
\end{equation}
The fact that the quantities on the right side of \eqref{eq:525} define probability measures on $\clv^d\times \RR_+$ follows from the identity $n^{-1} \int_0^{t_n} \psi_e(s) ds=1$.
From \eqref{eq:540} and chain rule for relative entropies (see \cite[Corollary 2.7]{buddupbook}), it follows that
\begin{equation}
	n^{-1} \sum\limits_{k=0}^{n-1} R\left(\bar \mu^{n,k+1} \| \rho(\bar L^{n,k+1})\right)
	= R\left(\beta^n \| \theta^n\right).
	\label{eq:541}
\end{equation}
With the identity in \eqref{eq:541}, the expectation on the right side of \eqref{eq:varrep} can be rewritten as
\begin{equation}
E\left[F(\bar L^n(t_n) ) + R\left(\beta^n \| \theta^n\right)\right].
\end{equation}
It is convenient to analyze the dynamics of $\bar L^n$ viewed backwards in time. Towards that end,  for each $n \in \NN$, define the $C(\RR_+: \clp(\Delta^o))$-valued random variable $\check \bfL^n$ by
\begin{equation} \label{eq:526}
\check \bfL^n(t) \doteq  \begin{cases} \bar L^n(t_n-t) & 0 \leq t \le t_n\\
\bar L^n(0) & t \ge t_n.\end{cases}
\end{equation}
Also, for each $n \in \NN$, define $\clm(\clv^d\times \RR_+)$-valued random variables $\check \Lambda^n$ by, for $A \subseteq \clv^d$ and $t \in \RR_+$,
\begin{equation} \label{eq:527}
	\check \Lambda^n(A \times [0,t]) \doteq \int_{t_n-t}^{t_n} \bar \Lambda^n(A\mid s) ds
	= \int_{0}^{t} \check\Lambda^n(A\mid s) ds,
\end{equation}
where $\bar \Lambda^n(A\mid s) \doteq 0$ for $s \le 0$, and $\check\Lambda^n(A\mid s) \doteq \bar \Lambda^n(A\mid t_n-s)$ for $s\in \RR_+$. 
For these time-reversed processes one can easily verify the following evolution equation:
for  $t \in \RR_+$ and $m(t) \le n$,
\begin{equation} \label{eq:timrev}
\begin{split}
\check \bfL^n (t) 
&= \check \bfL^n(0) - \int_{0}^{t} \sum\limits_{v\in\clv^d} v \check \Lambda^n(v|s)ds + \int_{t_n-t}^{t_n} \check \bfL^n(t_n - a(s)) ds.\\
\end{split}
\end{equation}
We now establish tightness of these time-reversed controlled processes and related collections of random variables.

\subsection{Tightness and Weak Convergence}\label{sec:prelimresults1}

We use the following estimate, for the difference between the harmonic series and the logarithm function, established in \cite{timtyr} (see also \cite{young} for a slightly  sharper estimate):
for any $n\ge 2$
\begin{equation}\label{eq:eulermaschestimate}
\gamma +\frac{1}{2(n+1)} < \sum\limits_{k=1}^n k^{-1} - \log n < \gamma + \frac{1}{2(n-1)},
\end{equation}
where $\gamma \approx 0.57721$ is the Euler-Mascheroni constant.

Recall the map $m: \RR_+ \to \NN_0$ defined by $m(t) \doteq \sup\{k \ge 1: t_k \leq t\}$. 
As an immediate consequence of the estimate in \eqref{eq:eulermaschestimate} and the observation that $t_n -s \leq t_{m(t_n-s)+1}$, we see that for all $n \in \NN$ and $s \in \RR_+$,
\begin{equation*}
\log(n+1 ) + \frac{1}{2(n+2)} - (s+1) \le t_{m(t_n-s)+1} - \gamma
\le \log(m(t_n-s)+2) + \frac{1}{2(m(t_n-s)+1)}.
\end{equation*}
The next lemma is a straightforward consequence of the above estimate (proof is omitted).
\begin{lemma}\label{lem:mtpsieasymptotics} 
For each $t \in \RR_+$, as $n \to \infty$, $n^{-1} m(t_n-t) \to \exp(-t).$
Additionally, for each $t \in \RR_+$, as $n \to \infty$,
\begin{equation}\label{eq:532to0}
\sup_{s\in [0,t]}\left| n^{-1} \psi_e(t_n-s) - \exp(-s) \right| \to 0.
\end{equation}

\end{lemma}

The next lemma establishes the tightness of the controlled processes introduced in Section \ref{sec:controlproc}.
\begin{lemma}\label{lem:lemtight}
The collection $\{(\check\bfL^n, \check \Lambda^n, \gamma^n, \beta^n, \theta^n), n \in \NN\}$
is tight in $C(\RR_+:\clp(\Delta^o))\times \clm(\clv^d\times \RR_+) \times (\clp(\clv^d \times\RR_+))^3$.
\end{lemma}

\begin{proof}
We begin by showing that $\{ \check \bfL^n , n \in \NN\}$ is tight. Since $\clp(\Delta^o)$ is compact, it suffices to show that there is a constant $C \in (0,\infty)$ such that for all $n \in \NN$ and $s,t \in \RR_+$,
$\| \check \bfL^n (t) - \check \bfL^n(s) \| \leq C |t - s|$.
Note that for all $k_1,k_2 \in \NN$ such that $k_1 < k_2 \leq t_n$, we have that $\|  \bar L^n(t_{k_1}) -  \bar L^n(t_{k_2})\| \leq 2 (t_{k_2} - t_{k_1})$. 
Using the linear interpolation property it then follows that for all $0 \le s \le t$,
$|\bar L^n(t)- \bar L^n(s)|\le 2 (t-s)$ which in turn shows that,
 for  $n \in \NN$ and $0\le s < t$,  $\| \check \bfL^n(t) - \check \bfL^n(s)\| \le 2|t-s|$.
Thus, it follows 
that $\{\check \bfL^n, \; n \in\NN\}$ is tight in $C(\RR_+: \clp(\Delta^o))$.

The tightness of $\{ \check \Lambda^n , n \in \NN\}$  in $\clm(\clv^d \times \RR_+)$
is immediate on observing that  for each $k \in \NN$  $\sup_{n \in\NN} \check \Lambda^n(\clv^d \times [0,k]) = k$.

Next, since $\clv^d$ is compact, the sequences $\{[\gamma^n]_1, \; n \in \NN\}$,  $\{[\beta^n]_1,\;  n \in \NN\}$, and  $\{[\theta^n]_1, \; n \in \NN\}$, are tight. Also, for each $n \in \NN$, $[\gamma^n]_2 = [\beta^n]_2 = [\theta^n]_2$, so it suffices to show that the sequence $\{ [\gamma^n]_2 ,\;  n \in \NN\}$ is tight. Observe that, for each $n \in \NN$,  if $n \ge m(t)$, then, since $[t_{m(t_n-t)+1}, t_n] \subseteq [t_n - t, t]$,
\begin{multline*}
[\gamma^n]_2([0,t]) = n^{-1} \int_0^t \psi_e(t_n-s)ds = n^{-1} \int_{t_n-t}^{t_n} \psi_e(s)ds\\
 \geq n^{-1} \sum\limits_{k=m(t_n-t)+1}^{n-1} \int_{t_k}^{t_{k+1}} \psi_e(s)ds = 1 -n^{-1}(m(t_n-t) +1).
\end{multline*}
From Lemma \ref{lem:mtpsieasymptotics},  for fixed $\veps>0$ and $t > \log (3 \veps^{-1})$,  we can find some  $n_0  > 3\veps^{-1}$ such that $n_0 \ge m(t)$ and
\[
\sup_{n\ge n_0}| n^{-1}m(t_{n}-t)- e^{-t}| \leq 3^{-1}\veps.
\]
Thus, 
\[
\inf\limits_{n\geq n_0}[\gamma^n]_2([0,t]) \geq 1- \veps.
\]
Since $\veps>0$ is arbitrary, the desired tightness follows.
\end{proof}

The next lemma gives a useful characterization for  the weak limit points of the tight collection in Lemma \ref{lem:lemtight}.

\begin{lemma}\label{lemchar}
	Let $(\check\bfL^*, \check \Lambda^*, \gamma^*, \beta^*, \theta^*)$
be a weak limit point of the sequence $(\check\bfL^n, \check \Lambda^n, \gamma^n, \beta^n, \theta^n)$.
Then, the following hold a.s.
\begin{enumerate}[(a)]
	\item The measure $\check \Lambda^*$ can be disintegrated as
	\[\check \Lambda^*(dv, ds) = \check \Lambda^*(dv\mid s) ds.
	\]
\item For all $t\in\RR_+$
	\begin{equation}\label{eq:eq617}
		\check \bfL^*(t) = \check \bfL^*(0) - \int_0^t \sum\limits_{v\in \clv^d} v \check \Lambda^*(v\mid s) ds + \int_0^t \check \bfL^*(s) ds.
	\end{equation}
	\item $\gamma^* = \beta^*$
\item For $t \in \RR_+$ and $x \in \Delta^o$,
\[
\gamma^*( \{e_x\}\times [0,t]) = \int_0^t \exp(-s) \check \Lambda^*(e_x\mid s) ds.
\]

\item For $t \in \RR_+$ and $x \in \Delta^o$,
\[
\theta^*( \{e_x\}\times [0,t]) = \int_0^t \exp(-s) \check \rho(\check \bfL^*(s))(e_x) ds.
\]

\end{enumerate}
\end{lemma}
\begin{proof}

\begin{enumerate}[(a)]
\item This is immediate on noting that for each $n \in \NN$ and $e_x \in \clv^d$, 
\[
\check \Lambda^n(e_x ,ds) = \check \Lambda^n(e_x \mid s) ds.
\]

\item Assume without loss of generality (by selecting the weakly convergent subsequence and appealing to Skorohod representation theorem)  that $\{(\check \bfL^n, \check \Lambda^n), \; n \in\NN\}$ converges almost surely to $(\check \bfL^*, \check \Lambda^*)$. For  $t \in \RR_+$ and $m(t) \le n$, recall the evolution equation \eqref{eq:timrev}.
Also note that
\begin{equation}\label{eq:weaklimitleqn1}
\int_{t_n-t}^{t_n} \check \bfL^n(t_n -a(s))ds
= \Bigg(\int_{t_n-t}^{t_n} \check \bfL^n(t_n -a(s))ds 
- \int_{t_n-t}^{t_n} \check \bfL^n(t_n -s)ds\Bigg)
+ \int_{0}^t \check \bfL^n(s)ds,
\end{equation}
and, for each $t \in \RR_+$,
\begin{equation}\label{eq:weaklimitleqn2}
\left\| \int_{t_n-t}^{t_n} \check \bfL^n(t_n -a(s))ds  - \int_{t_n-t}^{t_n} \check \bfL^n(t_n -s)ds\right\|
\le t \sup_{s \in[ t_n- t, t_n]} \| \bar L^n( a(s)) -  \bar L^n(s)\|.
\end{equation}
As in the proof of Lemma \ref{lem:lemtight},  for all $n \in \NN$ satisfying $t_n \geq t$ and $s \in [t_n-t,t_n]$, 
\begin{equation}\label{eq:weaklimitleqn3}
\| \bar{L}^n(a(s)) - \bar{L}^n(s)\| 
\leq  2 (m(t_n-t) +2)^{-1}.
\end{equation}
Combining \eqref{eq:weaklimitleqn1},  \eqref{eq:weaklimitleqn2}, and  \eqref{eq:weaklimitleqn3} and using the almost-sure convergence of $\{(\check \bfL^n, \check\Lambda^n), \; n\in\NN\}$ to $(\check \bfL^*, \check \Lambda^*)$ we see that, as $n \to \infty$, for each $t\in\RR_+$,
\[
\check \bfL^n(t) \to \check \bfL^*(0) - \int\limits_{0}^t \sum\limits_{v\in\clv^d} v \check \Lambda(v \mid s)ds + \int_0^t \check \bfL^*(s)ds,
\]
almost surely. The result follows.
\item
Fix $t >0$ and $x \in \Delta^o$. Then, for all $n\in \NN$ such that $t_n\ge t$,
\begin{equation}\label{eq:319}
\begin{split}
\gamma^n(\{e_x\} \times [0,t]) - \beta^n(\{e_x\} \times [0,t])
&= n^{-1} \int_0^{t} \psi_e(t_n-s) [ \bar \Lambda^n(e_x \mid t_n-s) -\bar \mu^n(e_x \mid t_n-s)] ds\\
&=  n^{-1} \int_{t_n-t}^{t_n} \psi_e(s) [ \bar \Lambda^n(e_x \mid s) -\bar \mu^n(e_x \mid s)] ds.
\end{split}
\end{equation}
Using the martingale-difference property, we see that for $1 \le l \le m$,
\begin{equation}\label{eq:squaredistmartdiff1}
E \left( n^{-1} \sum\limits_{k=l}^m [\delta_{\bar \nu^{n, k+1}}(e_x) - \bar \mu^{n, k+1}(e_x)]\right)^2
\le  \frac{m-l+1}{n^2}.
\end{equation}
Noting the identity
\begin{equation*}
n^{-1} \int_{t_l}^{t_{m+1}} \psi_e(s) [ \bar \Lambda^n(e_x \mid s) -\bar \mu^n(e_x \mid s)] ds
= n^{-1} \sum\limits_{k=l}^m [\delta_{\bar \nu^{n, k+1}}(e_x) - \bar \mu^{n, k+1}(e_x)],
\end{equation*}
we have from \eqref{eq:319} and \eqref{eq:squaredistmartdiff1} that for some $C_1 \in (0,\infty)$,
$$E[\gamma^n(\{e_x\} \times [0,t]) - \beta^n(\{e_x\} \times [0,t])]^2 \le C_1/n.$$
 The statement in (c) is now immediate.

\item \label{item:gammanlimit} As in part (b), without loss of generality, suppose that  
\[
\{(\check\bfL^n, \check \Lambda^n, \gamma^n, \beta^n, \theta^n), n \in \NN\},
\]
converges almost surely to $(\check\bfL^*, \check \Lambda^*, \gamma^*, \beta^*, \theta^*)$.
Fix $x \in \Delta^o$ and observe that for each $n \geq m(t)$,
\begin{multline}\label{eq:629triangle}
\left |n^{-1}\int_0^{t} \psi_e(t_n-s) \check \Lambda^n(e_x\mid s) ds - \int_0^{t} \exp(-s) \check \Lambda^*(e_x\mid s) ds\right |\\
\leq \left | n^{-1}\int_0^{t} \psi_e(t_n-s) \check \Lambda^n(e_x\mid s) ds -  \int_0^t \exp(-s)\check \Lambda^n(e_x \mid s)ds\right |\\ 
 + \left |  \int_0^t \exp(-s) \check \Lambda^n(e_x \times ds) -  \int_0^{t} \exp(-s)\check \Lambda^*(e_x\times ds) ds\right |.
\end{multline}
Next, 
\begin{multline*}
\left |n^{-1}\int_0^{t} \psi_e(t_n-s) \check \Lambda^n(e_x\mid s) ds -  \int_0^t \exp(-s) \check \Lambda^n(e_x \mid s)ds\right |\\
\leq \int_0^t \left| n^{-1}\psi_e(t_n-s) - \exp(-s)\right| \left | \check \Lambda^n (e_x \mid s)\right |ds
\leq  t \sup_{s\in[0,t]} \left| n^{-1} \psi_e(t_n-s) - \exp(-s)\right|,
\end{multline*}
and, by convergence of $\check \Lambda^n$ to $\check \Lambda^*$, as $n\to \infty$,
\begin{equation*}
\left |  \int_0^t \exp(-s) \check \Lambda^n(e_x \times ds) -  \int_0^{t}  \exp(-s) \check \Lambda^*(e_x\times ds) \right |
\to 0.
\end{equation*}
This, together with Lemma \ref{lem:mtpsieasymptotics} and \eqref{eq:629triangle},
shows that, as $n \to \infty$,
\[
n^{-1}\int_0^{t} \psi_e(t_n-s) \check \Lambda^n(e_x\mid s) ds \to \int_0^{t} \exp(-s) \check \Lambda^*(e_x\mid s) ds.
\]
On recalling the definition of $\gamma^n$ we now have the statement in (d).
\item Note that for each $t \in \RR_+$, $n \geq m(t)$, and $x \in \Delta^o$,
\begin{multline}\label{eq:thetanlimit668}
\left| \theta^n([0,t] \times \{e_x\}) - \int_0^t \exp(-s) \rho(\check \bfL^*(s)) (e_x)ds\right|\\
\leq  \Big| n^{-1} \int_0^t \psi_e(t_n-s) \rho(\bar{L}^n(a(t_n-s)))(e_x)ds
- \int_0^t \exp(-s)  \rho(\bar{L}^n(a(t_n-s)))(e_x)ds\Big|\\
+\Big|\int_0^t \exp(-s) \rho(\bar{L}^n(a(t_n-s)))(e_x)ds
- \int_0^t \exp(-s) \rho(\check \bfL^*(s)) (e_x)ds\Big|\\
\leq  \int_0^t \Big| n^{-1}\psi_e(t_n-s) - \exp(-s) \Big| ds
+ \int_0^t \exp(-s) \big|\rho(\bar{L}^n(a(t_n-s)))(e_x) - \rho(\check \bfL^*(s)) (e_x)\big|ds.
\end{multline}
From \eqref{eq:weaklimitleqn3},  for all $s \in [0,t]$, 
\begin{equation}\label{eq:thetanlimit676}
\big\|\rho(\bar{L}^n(a(t_n-s))) - \rho ( \check \bfL^*(s))\big\|
\le \big\|\bar{L}^n(a(t_n-s)) - \bar{L}^n(t_n-s) \big\|
\leq 2 (m(t_n-t)+2)^{-1} .
\end{equation}
Combining Lemma \ref{lem:mtpsieasymptotics}, \eqref{eq:thetanlimit668}, and \eqref{eq:thetanlimit676}, and sending $n\to \infty$, we now see that, for each $x \in \Delta^o$,  \[
\theta^*([0,t] \times \{e_x\}) = \int_0^t \exp(-s) \rho(\check \bfL^*(s)) (e_x)ds.
\]
The result follows.
\end{enumerate}
\end{proof}

\section{Laplace Upper Bound}\label{sec:laplaceupperbound}
The main result of the section is Theorem \ref{thm:laplaceupperbound1}, which gives the Laplace upper bound.
\begin{theorem}\label{thm:laplaceupperbound1}
	For every $F \in C_b(\clp(\Delta^o))$,
	\begin{equation*}
	\liminf_{n\to \infty}-n^{-1} \log E \exp[-n F(L^{n+1})] \ge \inf_{m \in \clp(\Delta^o)} [F(m) + I(m)].
	\end{equation*}
\end{theorem}

\begin{proof}
	Fix
	$F \in C_b(\clp(\Delta^o))$ and $\veps>0$. From the variational representation in  \eqref{eq:varrep}, for each $n \in \NN$ we can find $\{\bar \mu^{n,i}\} \in \Theta^n$ such that
	\begin{equation}\label{eq:lub1739}
		-n^{-1} \log E \exp[-n F(L^{n+1})]
		\ge E \left[F(\bar L^{n}(t_n)) +n^{-1} \sum\limits_{k=0}^{n-1} R\left(\bar \mu^{n,k+1}\big \| \rho(\bar L^{n,k+1})\right)\right] - \veps,
	\end{equation}
where the sequence $\{\bar L^{n, k}\}$ is defined by \eqref{eq:barlnk}.

For each $n \in \NN$, define the the $\clp(\Delta^o)$-valued continuous process $\bar L^n$ and random measures $\bar \Lambda^n, \bar \mu^n$ on $\clv^d \times [0, T_n]$ according to
\eqref{eq:522} and  \eqref{eq:519}, respectively. 
Also define, for each $n \in \NN$, $\gamma^n$, $\beta^n$, $\theta^n$, $\check \bfL^n$ and $\check \Lambda^n$ as in \eqref{eq:525},
\eqref{eq:526}, and \eqref{eq:527}, respectively. 
Recalling the identity in  \eqref{eq:541} we have that
\begin{equation}
		-n^{-1} \log E \exp[-n F(L^{n+1})] \ge E \left[F(\check \bfL^n(0)) +  R(\beta^n\| \theta^n)\right] - \veps.
\end{equation}
From Lemma \ref{lem:lemtight},
the collection $
\{ (\check\bfL^n, \check \Lambda^n, \gamma^n, \beta^n, \theta^n), \; n \in \NN\}$ is tight in $C(\RR_+:\clp(\Delta^o))\times \clm(\clv^d\times \RR_+) \times (\clp(\clv^d \times \RR_+))^3$.

Let $(\check\bfL^*, \check \Lambda^*, \gamma^*, \beta^*, \theta^*)$
be a weak limit point of the above sequence and suppose without loss of generality that the convergence holds along the full sequence and in the a.s. sense.
Disintegrating $\check \Lambda^*$   as in Lemma \ref{lemchar} (a), we see from part (b) of the same lemma
that, a.s., $\{\check \eta^*(\cdot \mid s), s \in \RR_+\} \in \clu(\check\bfL^*(0))$,
where 
\[
\check \eta^*(x \mid s) = \left(\sum\limits_{v \in \clv^d} v \check \Lambda^*(v\mid s)\right)_x.
\]
Also from parts (c), (d) and (e) of Lemma \ref{lemchar},
\begin{equation}\label{eq:615}
	 (\gamma^*, \theta^*) =  (\beta^*, \theta^*)
	 =  \left(\exp(-s) \check \Lambda^*(\cdot\mid s) ds, \exp(-s)  \rho(\check \bfL^*(s)) ds\right).
\end{equation}
Thus, using Fatou's lemma,
\begin{multline*}
		\veps + \liminf_{n\to \infty}-n^{-1} \log E \exp[-n F(L^{n+1})]
		\ge \liminf_{n\to \infty} E \left[F(\check \bfL^n(0)) +  R(\beta^n\| \theta^n)\right]\\
		\ge \left[F(\check \bfL^*(0)) +  R(\gamma^*\| \theta^*)\right]\
		= \left[F(\check \bfL^*(0)) +  R\left(\exp(-s) \check \Lambda^*(\cdot\mid s) ds\big\|  \exp(-s) \rho(\check \bfL^*(s)) ds\right)\right]\\
		 =  \left[F(\check \bfL^*(0)) +  \int_0^\infty \exp(-s) R\left(\check \Lambda^*(\cdot\mid s)
 		\big\|  \rho(\check \bfL^*(s)) \right) ds\right]\\
		\ge \left[F(\check \bfL^*(0)) + I(\check \bfL^*(0))\right]
		\ge \inf_{m \in \clp(\Delta^o)}\left[F(m) + I(m)\right],
\end{multline*}
where the second inequality  uses \eqref{eq:615} and the lower semicontinuity of relative entropy, the second identity uses the chain rule for relative entropies (see \cite[Corollary 2.7]{buddupbook}), and the last two inequalities use the fact that $\check \bfL^*(\cdot)$ solves \eqref{eq:eq617}, the relationship between $\check \eta^*$ and $\check \Lambda^*$, and the expression of the rate function $I$ given in \eqref{eq:equivrf}. The result follows on letting $\veps \to 0$.
\end{proof}

\section{Laplace Lower Bound}\label{sec:lowbd}
The main result of the section is Theorem \ref{thm:lowbd} which gives the  Laplace lower bound.
\begin{theorem}\label{thm:lowbd}
	For every $F \in C_b(\clp(\Delta^o))$,
	\begin{equation*}
	\limsup_{n\to \infty}-n^{-1} \log E \exp[-n F(L^{n+1})] \le \inf_{m \in \clp(\Delta^o)} [F(m) + I(m)].
	\end{equation*}
\end{theorem}

The proof of Theorem \ref{thm:lowbd} is postponed until Section \ref{sec:pflowbd}. In the next section we establish some estimates that are used in the proofs of Lemma \ref{lem:d1} and Lemma \ref{lem:d23}.

\subsection{Preliminary Estimates}\label{sec:prelimest}
Fix $F \in C_b(\clp(\Delta^o))$ and $\veps>0$.
In order to prove Theorem \ref{thm:lowbd} we can assume without loss of generality that $F$ is Lipschitz (see \cite[Corollary 1.10]{buddupbook}): for some $F_{\mbox{\tiny{lip}}} \in (0,\infty)$,
\[
|F(m) - F(\tilde m)| \le F_{\mbox{\tiny{lip}}} \|m - \tilde m\|, \;  m, \tilde m \in \clp(\Delta^o).
\]
Choose $m_0 \in \clp(\Delta^o)$ such that 
\begin{equation}
	F(m_0) + I(m_0) \le  \inf_{m \in \clp(\Delta^o)} [F(m) + I(m)] + \veps.
\end{equation}
Recalling the definition of the rate function, we choose $\eta_0 \in \clu(m_0)$ such that
\begin{equation}
	\int_0^{\infty} \exp(-s) R\left(\eta_0(\cdot \mid s) \| K(M_0(s)\right) ds
	\le I(m_0) + \veps, 
\end{equation}
where $M_0$ solves $\clu(m_0,\eta_0)$.

Now we make a series of approximations to $M_0$ and $\eta_0$
in order to get a more tractable and simple form near-optimal trajectory.

\subsubsection{Step 1: Ensuring Nondegeneracy.}
The relative entropy function $R(\cdot \| \theta)$ is not well behaved when a probability measure $\theta$ places small mass to points in its support. The following step addresses this problem.
Under Assumption \eqref{assu:kassumptions} there is a  unique stationary distribution for the transition probability matrix $A$, which we denote as $m_*$; note that $\inf\nolimits_{x\in\Delta^o}(m_*)_x>0$.
Let, for $s\in \RR_+$, $M_*(s) \doteq m_*$ and $\eta_*( \cdot  \mid s) \doteq  m_*$. Observe that $M_*$ solves $\clu(m_*,\eta_*)$. Define, for $\kappa \in (0,1)$ and $t \in \RR_+$,
\begin{align}\label{eq:mkappacomb}
	M_{\kappa}(t) \doteq (1-\kappa) M_0(t) + \kappa M_*(t),\nonumber\\
	\eta_{\kappa}(\cdot \mid t) \doteq
	(1-\kappa)\eta_{0}(\cdot \mid t) + \kappa \eta_{*}(\cdot \mid t),
\end{align}
and observe, with $m_{\kappa} \in \mathcal{P}(\Delta^o)$ defined as
$
m_{\kappa} \doteq (1-\kappa) m_0 + \kappa m_*,
$
that $M_{\kappa}$ solves $\clu(m_{\kappa},\eta_{\kappa})$.
Note also that
\begin{equation*}
\begin{split}
 \int_0^{\infty} \exp(-s) R\left(\eta_{\kappa}(\cdot \mid s) \| K(M_{\kappa}(s))\right) ds
  &\le (1-\kappa) \int_0^{\infty} \exp(-s) R\left(\eta_{0}(\cdot \mid s)  ds \| K(M_{0}(s) \right) ds\\
  &\le \int_0^{\infty} \exp(-s) R\left(\eta_{0}(\cdot \mid s)  ds \| K(M_{0}(s) \right) ds,
  \end{split}
  \end{equation*}
  where the first inequality follows from the convexity of relative entropy, the definitions of $m_{\kappa}, \eta_{\kappa}$, and $M_{\kappa}$ and the fact that
    \[
    \int_0^{\infty} \exp(-s) R\left(\eta_{*}(\cdot \mid s) \|  K(M_{*}(s))\right) ds=0.
    \]
  Let $\kappa_1 > 0$ be such that $  \| m_0 - m_{\kappa_1}\|  \le \min\{(F_{\mbox{\tiny{lip}}})^{-1} , 1\}\veps$, and, for  convenience, write $ (m_1, \eta_1,  M_1) \doteq (m_{\kappa_1}, \eta_{\kappa_1}, M_{\kappa_1})$.
 Then, 
  \begin{multline}
  F(m_{1}) + \int_0^{\infty} \exp(-s) R\left(\eta_{1}(\cdot \mid s) \| K(M_{1}(s))\right) ds\\
  \le F(m_0) + \int_0^{\infty} \exp(-s) R\left(\eta_{0}(\cdot \mid s)   \| K(M_{0}(s) \right) ds + \veps\\
  \le \inf_{m \in \clp(\Delta^o)} [F(m) + I(m)] + 3 \veps.\label{eq:358}
  \end{multline}
Also,  note that with $\delta_0 > 0$ defined as in Assumption \ref{assu:kassumptions},
\begin{equation}\label{eq:rhoboundlower}
K(m)(x) \ge \delta_0, \; x \in  \Delta^o,  m \in \clp(\Delta^o),
\end{equation}
which implies that for each $s \in \RR_+$,
\begin{equation}\label{eq:relentetarhom1}
	R\left(\eta_{1}(\cdot \mid s) \| K(M_{1}(s))\right) \le -\log \delta_0.
\end{equation}
Also, since $m_*$ is a positive vector, there is a $\delta>0$ such that for each $s \in \RR_+$ and $x \in \Delta^o$,
\begin{equation}\label{eq:1243}
	 M_{1}(s)(x) \ge \delta, \;  \eta_{1}(x \mid s) \ge \delta.
\end{equation}
Recall from the proof of the upper bound that the trajectories in the variational problem in the Laplace upper bound  are related to the controlled trajectories by  time reversal (see e.g., \eqref{eq:526}). Thus, we now introduce a time reversal of $M_1$, which, after further approximations, is then used to construct suitable controlled
trajectories.
Fix $T \in (0,\infty)$  large enough so that 
\begin{equation}\label{eq:sizeofT}
\exp(-T+1)|\log \delta_0| \le \veps
\end{equation}
Throughout the rest of the section, this $T$ is fixed.  Define, for $t \in [0,T]$,
\[
\hat M_1(t) \doteq M_1(T-t), \; \hat\eta_1(\cdot \mid t) \doteq \eta_1(\cdot \mid T-t),
\]
and note that, since $M_1$ solves $\clu(m_1, \eta_1)$, for $t \in [0,T]$,
\begin{equation}\label{eq:hatmit}
\begin{split}
M_1(T) + \int_0^t \hat\eta_1(s) ds - \int_0^t \hat M_{1}(s) ds 
&= \hat{M}_1(t).
 \end{split}
\end{equation}
Recalling the non-negativity of relative entropy, note that
\begin{equation}\label{eq:firscostapp}
\int_0^{\infty} \exp(-s) R\left(\eta_{1}(\cdot \mid s) \| K(M_{1}(s))\right) ds
\ge \exp(-T) \int_0^T \exp(s) R\left(\hat\eta_{1}(\cdot \mid s) \big\| K(\hat M_{1}(s))\right) ds.
\end{equation}
\subsubsection{Step 2: Continuity of Control}
Our next step mollifies the control $\hat \eta_1$ in a suitable manner so that it can be discretized at a later step.
For $\kappa >0$, define
\begin{equation}
	\hat\eta_1^{\kappa}(s) \doteq
	 \kappa^{-1}\int_{s}^{\kappa +s} \hat\eta_1(u) du, \; s \in [0,T],
\end{equation}
where $\hat \eta_1(u) \doteq \hat \eta_1(T)$ for $u\ge T$.
Also, define for $t \in [0,T]$,
\begin{equation}\label{eq:1153}
	\hat M_1^{\kappa}(t) \doteq M_1(T) + \int_0^t \hat\eta_1^{\kappa}(s) ds - \int_0^t \hat M_{1}^{\kappa}(s) ds.
\end{equation}
Note that there is a unique $\hat M_1^{\kappa} \in C([0,T]: \RR^d)$ that solves \eqref{eq:1153}, and that this $\hat M_1^{\kappa}$ satisfies, for each $s \in [0,T]$, 
\[
\sum\limits_{x \in \Delta^o}\hat M_1^{\kappa}(s)(x) = 1.
\]
We now show that for $\kappa$ sufficiently small, for each $ s \in [0,T]$ we have 
\[
\inf_{x\in\Delta^o}\hat M^{\kappa}_1(s)(x) > 0,
\]
namely that the solution to \eqref{eq:1153} in fact belongs to $C([0,T]: \clp(\Delta^o))$.
We can write, for $t \in [0,T]$,
\begin{equation}\label{eq:1134}
\hat M_1^{\kappa}(t) 
= M_1(T) + \int_0^{t}  \hat\eta_1(s)  ds 
- \int_0^t \hat M_{1}^{\kappa}(s) ds + \clr_1^{\kappa}(t),
\end{equation}
where 
\[
\clr_1^{\kappa}(t) \doteq \int_0^{t+\kappa}\hat\eta_1(u)\kappa^{-1}\int_{(u-k)^+}^udsdu - \int_0^t\hat\eta_1(u)du.
\]
Observe that, for each $t \in [0,T]$,
$\| \clr^{\kappa}_1(t)\| \le 3\kappa$. Combining this estimate with \eqref{eq:hatmit}
and \eqref{eq:1134}, we have,  for $t \in [0,T]$,
\begin{equation*}
\|  \hat{M}^{\kappa}_1(t) - \hat{M}_1(t)\| \le 3\kappa +  \int_0^t \|  \hat{M}^{\kappa}_1(s) - \hat{M}_1(s)\| ds,
\end{equation*}
from which we see, by an application of Gr\"{o}nwall's lemma, that
\begin{equation}\label{eq:mkappam1diffbound}
\sup_{0\le t \le T} \|\hat M_1^{\kappa}(t) - \hat M_1(t)\| \le 3\kappa \exp(T).
\end{equation}
We now assume that $\kappa$ is small enough so that, with $\delta$ as in \eqref{eq:1243}
\begin{equation}\label{eq:ledelt2}
3\kappa \exp(T) \le \delta/2.
\end{equation}
This, in view of \eqref{eq:1243}, ensures that $\hat M_1^{\kappa} \in C([0,T]: \clp(\Delta^o))$, and in fact
\begin{equation}\label{eq:945n}
\inf_{s \in [0,T], x \in \Delta^o}\hat M^{\kappa}_1(s)(x) \ge \delta /2. 
\end{equation}

Note that for $x\in \Delta^o$ and $u,s \in [0,T]$ satisfying $|u-s| \leq \kappa$,
\begin{equation}\label{eq:logmvtbound1}
| \log ( K(\hat M_1(u))(x)) - \log ( K(\hat M_1^{\kappa}(s))(x))|
\leq \delta_0^{-1}  |K(\hat M_1(u))(x))  - K(\hat M_1^{\kappa}(s))(x))|,
\end{equation}
so using the definition of relative entropy, \eqref{eq:mkappam1diffbound},  and  \eqref{eq:logmvtbound1} we see that
\begin{equation}\label{eq:eq218}
\begin{split}
R\left(\hat\eta_{1}(\cdot \mid u)) \big\| K(\hat M_{1}(u))\right) -  R\left(\hat \eta_{1}(\cdot \mid u)) \big\| K(\hat M_{1}^{\kappa}(s))\right)
&\leq \delta_0^{-1}\left(\| \hat{M}^{\kappa}_1(s) - \hat{M}_1(s) \| + \|\hat{M}_1(u) - \hat{M}_1(s)\|\right)\\
&\leq \delta_0^{-1} (3\kappa \exp(T) + 2\kappa).
\end{split}
\end{equation}
Using \eqref{eq:relentetarhom1} and the fact that $\bar \kappa^{-1} (1 - e^{-\bar \kappa}) \le  1$ for each $\bar \kappa \in (0,1)$, note that
\begin{multline}\label{eq:eq1010a}
\int_0^T \exp(s) \kappa^{-1} \int_s^{\kappa + s} R\left(\hat\eta_{1}(\cdot \mid u))  \big\| K(\hat M_{1}(u))\right)du ds\\
\le \kappa^{-1} \int_0^T R\left(\hat\eta_{1}(\cdot \mid u))  \big\| K(\hat M_{1}(u))\right)  (\exp(u) - \exp(u-\kappa))du\\
 + \kappa^{-1}  \int_T^{T+\kappa} |\log \delta_0| (\exp(u) - \exp(u-\kappa))du\\
\le \int_0^T \exp(u) R\left(\hat\eta_{1}(\cdot \mid u))  \big\| K(\hat M_{1}(u))\right)du  +  \kappa\exp({T} +{\kappa})|\log \delta_0|,
\end{multline}
so it follows from \eqref{eq:eq218}, \eqref{eq:eq1010a}, and convexity of relative entropy, that
\begin{multline}\label{eq:relentmkappam1011}
\exp({-T} )\int_0^T \exp(s) R\left(\hat\eta_{1}^{\kappa}(\cdot \mid s) \| K(\hat M_{1}^{\kappa}(s))\right) ds\\
\le \exp({-T})\int_0^{T} \exp(s) R\left(\hat \eta_{1}(\cdot \mid s)  \big\| K(\hat M_{1}(s))\right)  ds\\
+ \kappa\exp(\kappa)|\log\delta_0| + \delta_0^{-1}(3\kappa \exp(T) + 2\kappa).
\end{multline}
Now fix $\kappa_2$ small enough so that in addition to \eqref{eq:ledelt2} (with $\kappa = \kappa_2$) we have
\begin{equation*}
	\max\{ \kappa_2\exp(\kappa_2) | \log \delta_0| + \delta_0^{-1}(3\kappa_2 \exp(T) + 2\kappa_2),
	3\kappa_2 \exp(T)( 1 + F_{\mbox{\tiny{lip}}})\} \le \veps.
\end{equation*}
Henceforth, write $(\hat\eta_{2}, \hat M_{2}) \doteq (\hat\eta_{1}^{\kappa_2}, \hat M_{1}^{\kappa_2})$. Then, from
\eqref{eq:1153}, for $t \in [0,T]$,
\begin{equation}\label{eq:1233}
	\hat M_2(t) =  M_1(T) + \int_0^t \hat\eta_2(s) ds - \int_0^t \hat M_{2}(s) ds,
\end{equation}
and, from \eqref{eq:mkappam1diffbound} and our choice of $\kappa_2$, 
\begin{equation}
	\label{eq:ca}
\sup_{0\le t \le T} \|\hat M_2(t) - \hat M_1(t)\| \le \veps \min\{1, (F_{\mbox{\tiny{lip}}})^{-1}\}.
\end{equation}
Furthermore, from \eqref{eq:relentmkappam1011} and our choice of $\kappa_2$
\begin{multline}\label{eq:cb}
\exp(-T) \int_0^T \exp(s) R\left(\hat\eta_{2}(\cdot \mid s) \big\| K(\hat M_{2}(s))\right) ds	\\
\le  \exp(-T) \int_0^{T}\exp(s)   R\left(\hat \eta_{1}(\cdot \mid u)) \big\| K(\hat M_{1}(s))\right) ds + \veps.
\end{multline}
By construction, $\hat \eta_2$ is continuous and we can find  $C_1 \doteq C_1(\kappa_2) \in (0,\infty)$ such that
\begin{equation}\label{eq:220}
	\|\hat \eta_2(s) - \hat \eta_2(t)\|  \le C_1 |s-t|, \; \; s, t \in [0,T],
\end{equation}
and, with $\delta$ as in  \eqref{eq:1243}, for $s \in [0,T]$ and $x \in \Delta^o$,
\begin{equation} \label{eq:219}
	\hat\eta_{2}(x \mid s) \ge \delta.
\end{equation}
\subsubsection{Step 3: Piecewise Constant Approximation}
Now we carry out the last step in the approximation which is to replace continuous controls by piecewise constant controls.
For $\kappa >0$, define $\hat \eta_2^{\kappa}$ as
\begin{equation}
	\hat \eta_2^{\kappa} (\cdot \mid s) \doteq \hat \eta_2 (\cdot \mid j\kappa), \;\; s \in [j\kappa, (j+1)\kappa), \; j = 0, \ldots , \lfloor T \kappa^{-1} \rfloor.
\end{equation}
Also, for $t \in [0,T]$, let
\begin{equation}\label{eq:1208}
	\hat M_2^{\kappa}(t) =  M_1(T) + \int_0^t \hat\eta_2^{\kappa}(s) ds - \int_0^t \hat M_{2}^{\kappa}(s) ds,
\end{equation}
so that, with $\clr^{\kappa}_2(t) \doteq \int_0^t \hat\eta^{\kappa}_2(s)ds - \int_0^t \hat\eta_2(s)ds$,
we have that, for $t \in [0,T]$,
\begin{equation}\label{eq:1232}
\hat M_2^{\kappa}(t) =  M_1(T) + \int_0^t \hat\eta_2(s) ds - \int_0^t \hat M_{2}^{\kappa}(s) + \clr_2^{\kappa}(t).
\end{equation}
From \eqref{eq:220} and the definition of $\hat \eta^{\kappa}_2$, 
\[
\sup_{0\le t \le T}\clr_2^{\kappa}(t) \le T \sup_{s,t \in [0,T], |s-t| \le \kappa} \|\hat\eta_2(s) - \hat\eta_2(t)\|
\le C_1\kappa T,
\]
Combining the last estimate,  \eqref{eq:1233},  and  \eqref{eq:1232},  we have from Gr\"{o}nwall's lemma that
\begin{equation}\label{eq:221}
	\sup_{0\le t \le T} \|\hat M_2^{\kappa}(t) - \hat M_2(t)\| \le C_1\kappa T\exp({T}).
\end{equation}
Assume that $\kappa$ is sufficiently small so that $
C_1\kappa T\exp({T}) \le \delta/4.$
Then, from \eqref{eq:945n} it follows that $\hat M^{\kappa}_2$ is in $C([0,T]: \clp(\Delta^o))$.
Next, for $s \in [0,T]$,
\begin{multline*}
\left|R\left(\hat\eta_{2}(\cdot \mid s) \big\| K(\hat M_{2}(s))\right) -  R\left(\hat \eta_{2}^{\kappa}(\cdot \mid s) \big\| K(\hat M_{2}^{\kappa}(s))\right) \right|\\
\le 
\left|R\left(\hat\eta_{2}(\cdot \mid s) \big\| K(\hat M_{2}(s))\right) -  R\left(\hat \eta_{2}^{\kappa}(\cdot \mid s) \big\| K(\hat M_{2}(s))\right) \right| \\
 + \left|R\left(\hat\eta_{2}^{\kappa}(\cdot \mid s) \big\| K(\hat M_{2}(s))\right) -  R\left(\hat \eta_{2}^{\kappa}(\cdot \mid s) \big\| K(\hat M_{2}^{\kappa}(s))\right) \right|,
\end{multline*}
and, using \eqref{eq:rhoboundlower}, \eqref{eq:220}, and \eqref{eq:219}, 
\begin{multline}\label{eq:890a}
\left|R\left(\hat\eta_{2}(\cdot \mid s) \| K(\hat M_{2}(s))\right) -  R\left(\hat \eta_{2}^{\kappa}(\cdot \mid s) \| K(\hat M_{2}(s))\right) \right|\\
\le \sum\limits_{v\in\clv^d}\left( | \hat{\eta}_2(v \mid s) \log ( \hat \eta_2 (v \mid s) - \hat \eta^{\kappa}_2(v \mid s) \log( \hat \eta^{\kappa}_2 ( v \mid s)|
+ | \hat \eta_2(v \mid s) - \hat \eta^{\kappa}_s(v \mid s) | | \log ( K (\hat{M}_2(s)))|\right)\\
\le  \sum\limits_{v\in\clv^d} |  \hat \eta_2(v \mid s) - \hat \eta^{\kappa}_2(v \mid s)| ( | \log \delta_0 | + |\log \delta| + 1)
\le C_1 \kappa ( | \log \delta_0 | + |\log \delta| + 1).
\end{multline}
Using estimates analogous to those in  \eqref{eq:eq218} along with the estimate in  \eqref{eq:221}, we see that
\begin{equation*}
\left|R\left(\hat\eta_{2}^{\kappa}(\cdot \mid s) \big\| K(\hat M_{2}(s))\right) -  R\left(\hat \eta_{2}^{\kappa}(\cdot \mid s) \big\| K(\hat M_{2}^{\kappa}(s))\right) \right|
\le   \delta_0^{-1} \| \hat{M}^{\kappa}_2(s) - \hat M_2 (s)\| \le \delta_0^{-1} C_1 \kappa T \exp(T).
\end{equation*}
Combining the estimates in the last three displays we see that, for $s \in[ 0,T]$,
\begin{equation*}
\left|R\left(\hat\eta_{2}(\cdot \mid s)\big\| K(\hat M_{2}(s))\right) -  R\left(\hat \eta_{2}^{\kappa}(\cdot \mid s) \big\| K(\hat M_{2}^{\kappa}(s))\right) \right|
\le C_1 \kappa( |\log \delta_0| +|\log \delta| + 1 ) + \delta_0^{-1} C_1 \kappa T \exp(T).
\end{equation*}
It then follows that
\begin{multline}\label{eq:1206}
\exp(-T) \int_0^T \exp(s) R\left(\hat\eta_{2}^{\kappa}(\cdot \mid s)\big \| K(\hat M_{2}^{\kappa}(s))\right) ds\\
\le
\exp(-T)\int_0^T \exp(s) R\left(\hat\eta_{2}(\cdot \mid s) \big\| K(\hat M_{2}(s))\right) ds\\
+  C_1\kappa (|\log \delta_0| +|\log \delta|
+ 1) + \delta_0^{-1} C_1\kappa T \exp(T).
\end{multline}
Now, fix $\kappa_3 >0$ such that, in addition to $C_1 \kappa_3 T\exp(T) \le \delta/4$, 
\[
  C_1\kappa_{3}( (|\log \delta_0| +|\log \delta|
+ 1) + (\delta_0^{-1} + 1 + F_{\mbox{\tiny{lip}}})T \exp(T)) \le \veps.
\]
Writing,  $ (\hat\eta_{3}, \hat M_{3}) \doteq (\hat\eta_{2}^{\kappa_3}, \hat M_{2}^{\kappa_3})$,  we have from \eqref{eq:221}
\begin{equation}\label{eq:cc}
\sup_{0\le t \le T} \|\hat M_3(t) - \hat M_2(t)\| \le \veps \min\{1, (F_{\mbox{\tiny{lip}}})^{-1}\},
\end{equation}
and, from \eqref{eq:1208},
 for $t \in [0,T]$,
\begin{equation}\label{eq:1233b}
	\hat M_3(t) =  M_1(T) + \int_0^t \hat\eta_3(s) ds - \int_0^t \hat M_{3}(s) ds.
\end{equation}
Furthermore, from \eqref{eq:1206} and our choice of $\kappa_3$,
\begin{multline}\label{eq:cd}
\exp({-T}) \int_0^T \exp(s) R\left(\hat\eta_{3}(\cdot \mid s) \big\| K(\hat M_{3}(s))\right) ds\\
\le \exp(-T)\int_0^{T} \exp(s)   R\left(\hat \eta_{2}(\cdot \mid s)  \big\| K(\hat M_{2}(s))\right) ds + \veps.
\end{multline}
Finally, for each $j \le \lfloor T \kappa^{-1}\rfloor$, if $s \in [j\kappa, (j+1)\kappa)$, then $\hat\eta_3(s) = \hat\eta_3(\kappa j)$,
and from \eqref{eq:219},
\begin{equation} \label{eq:219b}
\inf_{x\in\Delta^o, s\in[0,T]}	\hat\eta_{3}(x \mid s) \ge \delta.
\end{equation}
Now we are ready to construct the control sequence.

\subsection{Constructing the Control Sequence}
 Let $T$ be large enough so that \eqref{eq:sizeofT} holds, and let
$q \doteq M_1(T)$, where $M_1$ is as in Section \ref{sec:prelimest}. From \eqref{eq:rhoboundlower}, we have that
\begin{equation}
	\sup_{m \in \clp(\Delta^o)} R(q \| K(m)) \le -\log \delta_0.
\end{equation}
Let $\veps_0>0$ be such that 
\begin{equation}\label{eq:vepsbound1}
\max\{1, F_{\mbox{\tiny{lip}}}\} \exp(T) \veps_0 ( 1 + \delta_0 T)
\le \veps.
\end{equation}
 Let $\{Y_i, \; i\in \NN\}$ be iid $\Delta^o$-valued random variables with common distribution $q$.  Using the law of large numbers we can find  $n_0 \in \NN$ such that 
\begin{equation}\label{eq:515}
	P\left(\left\|n^{-1} \left( \delta_{x_0} + \sum\limits_{i=1}^{n-1} \delta_{Y_i}\right) - q\right\| \ge \veps_0\right) \le \veps, \; n \ge n_0.
\end{equation}
\begin{construction}\label{const:controlseq1}
	Fix $n  \ge m(T)$, and let $c \doteq \kappa_3$, where $\kappa_3$ is defined as in Section \ref{sec:prelimest}.
	For $s \in \RR_+$, define $\hat \Lambda_3(s) \in \clp(\clv^d)$ as
	\[
	\hat \Lambda_3(s)(e_x) \doteq \hat \eta_3(s)(x), \; x \in \Delta^o.
	\]
\begin{enumerate}[(i)]
	\item Let $a_0(n) \doteq m(t_n-T)$, so that $t_n-T \in [t_{a_0(n)}, t_{a_0(n)+1})$. 
	Assume that $n$ is large enough so that $a_0(n) \ge n_0$.	
	For $1 \le i \le a_0(n)+1$, let $\bar \nu^{n,i} \doteq \delta_{Y_i}$ and $\bar \mu^{n,i} \doteq q$, define the random measures $\bar L^{n,i}$ recursively as in \eqref{eq:barlnk} using this choice of $\bar \nu ^{n,i}$.
	Also, let $\bar \clf^{n,i} \doteq \sigma\{\bar L^{n,j}, \; 1 \le j \le i\}$.
	\item Consider the set 
	\begin{equation*}A_n \doteq \left\{\{\|\bar L^{n,a_0(n)+2}  - q\| \ge \veps_0\right\}
	= \left\{\left\|(a_0(n)+2)^{-1} \left( \delta_{x_0} + \sum\limits_{i=1}^{a_0(n)+1} \delta_{Y_i}\right) - q\right\| \ge \veps_0\right\}.
	\end{equation*}
	For $i \ge a_0(n)+2$, define the random measure $\bar \mu^{n,i} $, together with $\bar \nu^{n,i}$ and $\bar L^{n,i}$, recursively as in \eqref{eq:barlnk}, as follows.
	
	\begin{itemize}
	 \item On  $A_n \in \bar \clf^{n,a_0(n)+2}$, let, for $i \ge a_0(n) + 2$,  $\bar \mu^{n,i} \doteq \rho(\bar L^{n,i})$.
	 \item On  $A_n^c$,  define $\bar \mu^{n,i}$, for $i \ge a_0(n) + 2$, as follows. Let $\sigma_n^+ \doteq t_{a_0(n)+1}$ and for $j =0, 1, \ldots, \lfloor Tc^{-1}\rfloor$ and $t_i \in [\sigma_n^+ + jc, \sigma_n^++ (j+1)c)$, define
	$\bar \mu^{n,i} \doteq \hat \Lambda_3(cj).$
	For $i \ge m(\sigma_n^+ + (\lfloor Tc^{-1}\rfloor+1)c)$, define $\bar \mu^{n,i} \doteq \rho(\bar L^{n,i}).$
	\end{itemize}
	\item Using the above $\{\bar L^{n,i}\}$, define $\bar L^n(\cdot)$ using \eqref{eq:411} and  random measures $\bar\Lambda^n$ and $\bar \mu^n$ on $\clv^d \times [0,t_n]$ using
	\eqref{eq:519}.
	\item For each $n \in \NN$, define the continuous process $\check \bfL^n(\cdot)$ using \eqref{eq:526} and the  $\clm(\clv^d\times \RR_+)$-valued random variable $\check \Lambda^n$ using
\eqref{eq:527}.
\item For each $n \in \NN$, define  $\clp(\clv^d \times \RR_+)$-valued random variables $\gamma^n, \beta^n, \theta^n$ using
\eqref{eq:525}.
\end{enumerate}
\end{construction}

With $\{\bar \mu^{n,k}\}$ and $\{\bar L^{n,k}\}$ sequences defined as above, we have from the variational representation in \eqref{eq:varrep}
\begin{equation}\label{eq:varreplow}
	-n^{-1}\log E \exp[-n F(L^{n+1})]\\
	\le  E \left[F(\bar L^n(t_n) ) + n^{-1} \sum\limits_{k=0}^{n-1} R\left(\bar \mu^{n,k+1} \| \rho(\bar L^{n,k+1})\right)\right].
\end{equation}

We begin with the following lemma. 
\begin{lemma}
	\label{lem:d1}
	We have that, with $d_1 = 2\|F\|_{\infty}+1$,
\[ \limsup_{n\to \infty} E(F(\check \bfL^n(0))) \le F(\hat M_3(T)) + d_1 \veps.
\]
\end{lemma}
\begin{proof}
Recall $\sigma_n^+\doteq t_{a_0(n)+1}$. Also
define sequences $\{\sigma_n\}$, $\{\sigma_n^-\}$ by
\[
\sigma_n \doteq t_n-T, \;\; \sigma_n^- \doteq t_{a_0(n)}, \;\; n \in \NN,
\]
and observe that $\sigma_n \in [t_{a_0(n)}, t_{a_0(n)+1})$. For each $n \in \NN$ and $s \in \RR_+$, define also  $\sigma_{n,s} \doteq s + \sigma_n$. For each $n \in \NN$ and $t \in [0,T]$, define
\begin{equation}\label{lem:mbardeflem7}
\bar M^n(t) \doteq \bar L^n(\sigma_{n,t}),
\end{equation}
 and observe that $\check \bfL^n(0) = \bar L^n(t_n) = \bar M^n(T)$.
For each $t \in [0,T]$ and $x \in \Delta^o$, let $\bar \eta^n(t)(x) \doteq \bar \Lambda^n(e_x\mid \sigma_{n,t})$ and
\[
\clr^n(t) \doteq \int_0^t \bar{M}^n(s)ds - \int_0^t \bar L^n(a(\sigma_{n,s}))ds.
\]
Note from \eqref{eq:522} that, for each $t \in [0,T]$,
\begin{equation}
	\bar M^n(t) 
	= \bar M^n(0) + \int_0^t \bar \eta^n(s) ds - \int_0^t \bar  M^n(s) ds + \clr^n(t). \label{eq:829}
\end{equation}
From an estimate as in  the proof of Lemma \ref{lem:lemtight}, we see that
\begin{multline*}
\sup_{0\leq t \leq T} \| \clr^n(t)\| 
=  \sup_{0\leq t \leq T} \left\| \int_0^t \left(\bar L^n(\sigma_{n,s}) - \bar L^n(a(\sigma_{n,s}))\right)ds\right\|
\le 2    \int_0^T |\sigma_n+s - a(\sigma_{n,s})| ds \\
\le 2 \sum\limits_{j= a_0(n)}^{n-1} \int_{t_j}^{t_{j+1}} | s - t_j|ds 
 \le 2 a_0(n)^{-1}.
\end{multline*}
Let $\clr^n_1(t) \doteq \int_0^{\sigma_n^+ - \sigma_n} \bar \eta^n(s) ds + (\bar M^n(0) - \bar L^n(\sigma_n^-)) + \clr^n(t)$,
so that
\begin{equation}\label{eq:445}
\bar M^n(t)
= \bar L^n(\sigma_n^-) + \int_{\sigma_n^+-\sigma_n}^t \bar \eta^n(s) ds - \int_0^t \bar  M^n(s) ds + \clr^n_1(t).
\end{equation}
Note that, as in the proof of Lemma \ref{lem:lemtight},
 \begin{equation}\label{eq:clr1bound}
\sup_{0\leq t \leq T}\| \clr^n_1(t)\| \leq (\sigma_n^+ - \sigma_n) + \| \bar L^n(\sigma_n)  - \bar L^n(\sigma_n^-)\| +  2 a_0(n)^{-1}
\le  5a_0(n)^{-1}.
\end{equation}
For $n \in \NN$ and $t \in [0,T]$,  define
\[
\clr^n_2(t) \doteq \int_{\sigma_n^+-\sigma_n} ^t\bar \eta^n(s) ds -  \int_{0} ^t\hat \eta_3(s) ds.
\]
Then,  for $t \in [0,T]$,
\begin{equation}\label{eq:450}
\int_{\sigma_n^+-\sigma_n} ^t\bar \eta^n(s) ds
= \int_{0} ^t\hat \eta_3(s) ds + \clr^n_2(t).
\end{equation}
For $t \in [0,T]$, define 
\begin{align*}
\clr^{n,1}_2(t) &\doteq \int_{\sigma_n^+-\sigma_n}^t ( \bar \Lambda^n(\sigma_{n,s}) - \bar \mu^n(\cdot \mid \sigma_{n,s})) ds\\  \clr^{n,2}_2(t) &\doteq \int_{\sigma_n^+-\sigma_n}^t ( \bar\mu^n(\cdot \mid \sigma_{n,s}) - \hat \Lambda_3(s))ds,
\end{align*}
and let $\clr^{n,3}_2 \doteq - \int_0^{\sigma_n^+ - \sigma_n}\hat\eta_3(s)ds$, so that
\[
\clr^n_2 (t) = \clr^{n,1}_2(t) + \clr^{n,2}_2(t) + \clr^{n,3}_2.
\]
Observe that for each $t \in [0,T]$, on $A_n^c$
\begin{align}\label{eq:remainderoverlapintervalsdiff}
\| \clr^{n,2}_2(t)\| &\le 2  (\lfloor T c^{-1} \rfloor + 1)a_0(n)^{-1} \nonumber\\
\| \clr^{n,3}_2\|& \le  (\sigma_n^+ - \sigma_n) \le  a_0(n)^{-1},
\end{align}
where the first inequality follows on recalling the definition of $\bar \mu^n$ in Construction \ref{const:controlseq1}.
Now, consider random signed measures $\{\bar \Delta^{n,i} , n \in \NN, i \in\NN_0\}$ and $\{ \bar\Delta^n(s), n \in \NN, s \in \RR_+\}$ defined by
\begin{align*}
\bar \Delta^{n,i} &\doteq \delta_{\bar \nu^{n,i}} - \bar \mu^{n,i}, \quad n \in \NN, i \in \NN_0\\
\bar\Delta^n(s) &\doteq \bar \Lambda^n(s) - \bar \mu(\cdot \mid s), \quad n \in \NN, s \in \RR_+,
\end{align*}
and observe that, for each $n \in \NN$ and $x \in \Delta^o$, $\{ \bar \Delta^{n,i}(e_x)\}_{i\in\NN_0}$ is a martingale difference sequence. Thus,  using Burkholder's inequality, we see that, for all $n \in \NN$,
\begin{multline*}
E\left( \sup_{0\leq t \leq T} \| \clr^{n,1}_2(t) \|^2\right) =  E\left(\sup_{0\leq t \leq T} \left\| \int_{\sigma_n^+}^{t+\sigma_n} \bar \Delta^n(s)ds\right\|^2 \right)
\le 4 d \sum\limits_{i=a_0(n)+1}^{m(\sigma_{n,T})+1}E\left[(i+2)^{-2} \|  \delta_{\bar \nu^{n,i+1}} - \bar \mu^{n,i+1}\|^2\right]\\
\le 4d \sum\limits_{i=a_0(n)+1}^{m(\sigma_{n,T})+1} (i+2)^{-2} \le 4d (a_0(n))^{-1}.
\end{multline*}
From \eqref{eq:remainderoverlapintervalsdiff} it then follows that
\begin{equation}\label{eq:clr2bound}
\sup_{0\le t \le T} \|\clr^n_2(t)\| \le  a_0(n)^{-1}(2\lfloor Tc^{-1}\rfloor +4) + r_n,
\end{equation}
where $r_n\to 0$ in $\mathcal{L}^2$ as $n \to \infty$.
Thus, with $\clr^n_3(t) \doteq ( \bar L^n(\sigma_n^-) - q)  + \clr^n_1(t) + \clr^n_2(t)$,
we have from \eqref{eq:445} and \eqref{eq:450} that, for $t \in \RR_+$,
\begin{align*}
\bar M^n(t) = q + \int_{0} ^t\hat \eta_3(s) ds   - \int_0^t \bar  M^n(s) ds + \clr^n_3(t).
\end{align*}
From \eqref{eq:clr1bound} and \eqref{eq:clr2bound}, on the set $A_n^c$,
\begin{equation*}
\begin{split}
\sup_{0\leq t\leq T} \| \clr^n_3(t)\| &\leq a_0(n)^{-1} \left(5+  (2\lfloor Tc^{-1}\rfloor +4)\right)
+ r_n+ \|q - \bar{L}^{n, a_0(n)+1}\|\\
&\le  a_0(n)^{-1}(9 + 2\lfloor Tc^{-1}\rfloor)+ r_n  + \veps_0.
\end{split}
\end{equation*}
For notational convenience write $\tilde{C}_{T,c} \doteq 9 + 2\lfloor T c^{-1} \rfloor$,
so that
\[
\sup_{0\le t \le T} \| \clr^n_3(t)\| \le a_0(n)^{-1} \tilde{C}_{T,c} + r_n + \veps_0.
\]
Now, from \eqref{eq:1233b}, Gr\"{o}nwall's lemma, and the fact that $M_1(T) = q$,   on $A_n^c$,
\begin{equation}\label{eq:516b}
	\sup_{0\le t \le T}\|\bar M^n(t) - \hat M_3(t)\|
	\le \exp(T) \left(a_0(n)^{-1}\tilde{C}_{T,c}+  r_n + \veps_0\right).
\end{equation}
Finally,
\begin{multline*}
\limsup_{n\to \infty} E(F(\check \bfL^n(0))) \le  \limsup_{n\to \infty} E(F(\check \bfL^n(0))1_{A_n^c}) + \|F\|_{\infty} \veps
= \limsup_{n\to \infty} E(F(\bar M^n(T))1_{A_n^c}) + \|F\|_{\infty} \veps\\
\le  \limsup_{n\to \infty} E(F(\hat M_3(T))1_{A_n^c})
 + \|F\|_{\infty} \veps +  \limsup_{n\to \infty} F_{\mbox{\tiny{lip}}}\exp(T)\left(a_0(n)^{-1}\tilde{C}_{T,c}+ Er_n+ \veps_0\right)\\
\le F(\hat M_3(T)) + 2\|F\|_{\infty} \veps + F_{\mbox{\tiny{lip}}}\exp(T)\veps_0
\le  F(\hat M_3(T)) + (2\|F\|_{\infty}+1) \veps,
\end{multline*}
where the first inequality uses \eqref{eq:515}, the identity uses the observation below \eqref{lem:mbardeflem7}, the second inequality uses \eqref{eq:516b} and the last inequality
 is due to \eqref{eq:vepsbound1}.
\end{proof}

The next lemma estimates the cost of our constructed controls.
\begin{lemma}\label{lem:d23} 
Let $\{\beta^n, \theta^n, \; n \in \NN\}$ be as in Construction \ref{const:controlseq1}.
Then,
	 \begin{equation*}
	 \limsup_{n\to \infty} E\left[ R(\beta^n\| \theta^n)\right]\\
	 \le 
	\exp(-T) E \int_0^{T}  \exp(s) R\left(\hat \eta_{3}(\cdot \mid s) \| K(\hat M_{3}(s))\right) ds + 2\veps.
	\end{equation*}
\end{lemma}
\begin{proof}
	From \eqref{eq:541}, for each $n \in \NN$,
	\begin{equation}\label{eq:522b}
		R(\beta^n\| \theta^n)
		= n^{-1} \sum\limits_{k=0}^{a_0(n)} R\left(\bar \mu^{n,k+1} \| \rho(\bar L^{n,k+1})\right)
		+ n^{-1} \sum\limits_{k=a_0(n)+1}^{n-1} R\left(\bar \mu^{n,k+1} \| \rho(\bar L^{n,k+1})\right).
	\end{equation}
Recall from our choice of $T$ made below \eqref{eq:1243} that
$\exp(-T+1) |\log \delta_0| \le \veps$, so, from Lemma \ref{lem:mtpsieasymptotics}, we have that for all sufficiently large $n$,
\begin{equation*}
n^{-1} \sum\limits_{k=0}^{a_0(n)} R\left(\bar \mu^{n,k+1}  \big\| \rho(\bar L^{n,k+1})\right)
\le n^{-1}(a_0(n)+1) |\log \delta_0| 
\le \exp(-T+1) |\log \delta_0| \le \veps.
\end{equation*}
Next, note that, by Construction \ref{const:controlseq1}, on the set $A_n$,
\[
n^{-1} \sum\limits_{k=a_0(n)+1}^{n-1} R\left(\bar \mu^{n,k+1} \big\| \rho(\bar L^{n,k+1})\right)=0.
\]
Now we estimate the above relative entropy on the set $A_n^c$.
Along the lines of \eqref{eq:540},
\begin{equation*}
n^{-1} \sum\limits_{k=a_0(n)+1}^{n-1} R\left(\bar \mu^{n,k+1}\big \| \rho(\bar L^{n,k+1})\right)
	 n^{-1} \int_{\sigma_n^+-\sigma_n}^{T}
	\psi_e(\sigma_{n,s}) R\left(\bar \mu^n(\cdot \mid \sigma_{n,s}) \big\|  \rho(\bar L^n(a(\sigma_{n,s})))\right) ds.
\end{equation*}
Let 
\begin{equation*}
A_n(s) \doteq \big| R\left( \bar \mu^n(\cdot \mid  \sigma_{n,s}) \big\| \rho (\bar L^n(a( \sigma_{n,s})))\right)
- R\left(\hat \Lambda_3(\cdot  \mid s ) \big\| \rho (\bar{L} ^n(a(\sigma_{n,s})))\right)\big|,
\end{equation*}
and
\[
B_n(s) \doteq \big| R\left(\hat \Lambda_3(\cdot \mid s )\big \| \rho (\bar{L} ^n(a(\sigma_{n,s})))\right) - R \left(\hat \Lambda_3(\cdot \mid s) \big\| \rho(\bar L^n( \sigma_{n,s}))\right)\big|.
\]
Then, recalling the definition of $\bar{M}^n$ from \eqref{lem:mbardeflem7}, we have
\begin{multline}
n^{-1} \int_{\sigma_n^+-\sigma_n}^{T}
	\psi_e(\sigma_{n,s}) R\left(\bar \mu^n(\cdot \mid \sigma_{n,s}) \|  \rho(\bar L^n(a(\sigma_{n,s})))\right) ds \\
	 \le n^{-1} \int_{\sigma_n^+-\sigma_n}^{T}
	\psi_e(\sigma_{n,s}) R\left( \hat \Lambda_3(\cdot \mid s) \| \rho(\bar M^n(s))\right) ds  
	+n^{-1} \int_{\sigma_n^+-\sigma_n}^{T}
	\psi_e(\sigma_{n,s}) [A_n(s) + B_n(s) ]ds.
\end{multline}
Estimates analogous to those in \eqref{eq:890a} show that
\[
A_n(s) \le  \left( | \log \delta_0| + | \log \delta| + 1\right)  \| \bar \mu^n( \cdot \mid \sigma_{n,s}) - \hat \Lambda_3(\cdot \mid s)\| , 
\]
so it follows, on recalling the definition of $\bar \mu^{n,i}$ from
Construction \eqref{const:controlseq1}, that, on $A_n^c$,
\begin{equation*}
\begin{split}
n^{-1} \int_{\sigma_n^+-\sigma_n}^{T}
	\psi_e(\sigma_{n,s} )A_n(s) ds
	&\le   \left( | \log \delta_0| + | \log \delta| + 1\right) \int_{\sigma_n^+ - \sigma_n}^{T} \| \bar \mu^{n}( \cdot \mid  \sigma_{n,s}) - \hat \Lambda_3(\cdot \mid s) \|ds\\
	&\le  2a_0(n)^{-1}  \left( | \log \delta_0| + | \log \delta| + 1\right) ( \lfloor T c^{-1} \rfloor + 1) ,
	\end{split}
\end{equation*}
where the final inequality follows as the estimate in \eqref{eq:remainderoverlapintervalsdiff}. Next, using estimates similar to \eqref{eq:weaklimitleqn3} and \eqref{eq:eq218}, we  see that 
for $n\in\NN$ and  $s \in [\sigma^+_n- \sigma_n, T]$,
\[
B_n(s) \le \delta_0^{-1} \| \bar L^n(a(\sigma_{n,s})) - \bar L^n(  \sigma_{n,s})\| \le \delta_0^{-1} 2 a_0(n)^{-1}.
\]

It now follows from Lemma \ref{lem:mtpsieasymptotics}  that as $n \to \infty$,
\[
n^{-1} \int_{\sigma^+_n-\sigma_n}^T \psi_e(\sigma_{n,s})B_n(s)ds \to 0.
\]
From an estimate analogous to \eqref{eq:eq218}, and   \eqref{eq:vepsbound1} and \eqref{eq:516b}, it follows that, 
\begin{multline*}
n^{-1} \int_{\sigma_n^+-\sigma_n}^{T}
	\psi_e(\sigma_{n,s}) R\left(\hat \Lambda_3(\cdot \mid s) \big\|  \rho(\bar M^n(s))\right) ds\\
	\le n^{-1} \int_{\sigma_n^+-\sigma_n}^{T}
	\psi_e(\sigma_{n,s}) R\left(\hat \Lambda_3(\cdot \mid s) \big\|  \rho(\hat M_3(s))\right) ds + \veps + \clr(n),
\end{multline*}
where $\clr(n) \to 0$ in $\mathcal{L}^2$ as $n\to \infty$.
Also, from Lemma \ref{lem:mtpsieasymptotics}, it follows that \[
\sup_{0\le s \le T} |n^{-1} \psi_e(\sigma_{n,s}) - \exp({-(T-s)})| \to 0 \mbox{ as }
\]
as $n\to \infty$.
Combining the preceding bounds and convergence, we have, on $A_n^c$, for sufficiently large $n$,
\begin{equation*}
n^{-1} \sum\limits_{k=a_0(n)+1}^{n-1} R\left(\bar \mu^{n,k+1} \big\| \rho(\bar L^{n,k+1})\right)
\le \exp(-T) \int_0^T \exp(s) R\left(\hat \Lambda_3(\cdot \mid s) \|  \rho(\hat M_3(s))\right) ds + \veps + \clr_1(n),
\end{equation*}
where $\clr_1(n) \to 0$ in $\mathcal{L}^2$ as $n\to \infty$.
Finally, taking expectations in \eqref{eq:522b},
\begin{equation*}
 E\left[ R(\beta^n\| \theta^n)\right]
 \le \exp({-T})\int_0^T \exp(s) R\left(\hat \Lambda_3(\cdot \mid s) \|  \rho(\hat M_3(s))\right) ds + 2 \veps + E\clr_1(n),
\end{equation*}
and the result follows on sending $n\to \infty$, and recalling the relations between $\hat \eta_3$ and $\hat \Lambda_3$, and $\rho(m)$ and $K(m)$, $m \in \clp(\Delta^o)$.
\end{proof}

We now complete the proof of the Laplace lower bound.

\subsection{\bf Proof of Theorem \ref{thm:lowbd}.}\label{sec:pflowbd}
From \eqref{eq:358} and \eqref{eq:firscostapp}, and recalling that $\hat M_1(T)= m_1$,
  \begin{equation}
  F(\hat M_{1}(T)) + \exp({-T})\int_0^{T} \exp(s) R\left(\hat \eta_{1}(\cdot \mid s) \big\| K(\hat M_{1}(s))\right) ds
  \le \inf_{m \in \clp(\Delta^o)} [F(m) + I(m)] + 3 \veps.\label{eq:359}
  \end{equation}
  Combining this with \eqref{eq:ca}, \eqref{eq:cb}, \eqref{eq:cc}, \eqref{eq:cd} we have
  \begin{equation}
  F(\hat M_{3}(T)) + \exp({-T})\int_0^{T} \exp(s) R\left(\hat \eta_{3}(\cdot \mid s) \big\| K(\hat M_{3}(s))\right) ds
  \le \inf_{m \in \clp(\Delta^o)} [F(m) + I(m)] + 7 \veps.\label{eq:402}
  \end{equation}
  Using Lemma \ref{lem:d1}, the identity in \eqref{eq:541}, and the inequality in \eqref{eq:varreplow},
  \begin{equation*}
  \begin{split}
	\limsup_{n\to \infty}-n^{-1} \log E \exp[-n F(L^{n+1})]
	&\le \limsup_{n\to \infty} E \left[F(\check \bfL^n(0)) +  R(\beta^n\| \theta^n)\right]\\
		&	\le F(\hat M_3(T)) +  \limsup_{n\to \infty} E\left[R(\beta^n\| \theta^n)\right] + d_1\veps.
	\end{split}
  \end{equation*}
  Also, from Lemma \ref{lem:d23},
  \begin{equation}\label{eq:u1}
  	\limsup_{n\to \infty} E\left[R(\beta^n\| \theta^n)\right] 
	\le \exp({-T})\int_0^T \exp(s) R\left(\check\eta_{3}(\cdot \mid s) \big\| K(\check M_{3}(s))\right) ds + 2\veps.
	\end{equation}
	Combining the above bounds we have that
	\begin{equation*}
		\limsup_{n\to \infty}-n^{-1} \log E \exp[-n F(L^{n+1})] 
		\le \inf_{m \in \clp(\Delta^o)} [F(m) + I(m)] + 9 \veps + d_1\veps.
	\end{equation*}
Since $\veps>0$ is arbitrary, the result follows.
  \hfill \qed
 \section{Compactness of level sets.}\label{sec:levelset}
 In this section we show that the function $I$ defined in \eqref{def:ratefunction1} is a rate function. For this it suffices to show that for any $k \in (0,\infty)$, the set
 $S_k = \{m \in \clp(\Delta^o): I(m)\le k\}$ is compact in $\clp(\Delta^o)$.  Let $\{m_n, \; n \in \NN\}$ be a sequence in $S_k$. Since $\clp(\Delta^o)$ is compact, $\{m_n\}$ converges along a subsequence to some limit point $m \in \clp(\Delta^o)$. It suffices to show that $m \in S_k$. Since $m_n \in S_k$, for each $n \in \NN$ we can find $\eta^n \in \clu(m_n)$
 such that
 \begin{equation}\label{eq:935}
 \int_0^{\infty} \exp(-s) R\left(\eta^n(\cdot \mid s) \| K(M^n(s)\right) ds
 \le I(m_n) +n^{-1}  \le k +n^{-1},
 \end{equation}
 where $M^n$ solves $\clu(m^n, \eta^n)$.
 Define probability measures $\hat \gamma^n, \hat \theta^n$ on
 $\Delta^o \times \RR_+$ as, for $t \in \RR_+$ and $x \in \Delta^o$,
 \begin{align*}
 \hat \gamma^n(\{x\}\times [0,t]) &= \int_0^t e^{-s} \eta^n(x \mid s) ds\\
 \hat \theta^n(\{x\}\times [0,t]) &= \int_0^t e^{-s} K(M^n(s))(x) ds.
 \end{align*}
 Since $\Delta^o$ is compact and $[\hat \gamma^n]_2(ds) = [\hat \theta^n]_2(ds) = e^{-s} ds$, it follows that the sequences $\{\hat \gamma^n, \; n\in \NN\}$, $\{\hat \theta^n, n \in \NN\}$ are tight in 
 $\clp(\Delta^o \times \RR_+)$. Consider a further subsequence (of the subsequence along which $m^n$ converges) along which $\hat \gamma^n$
 and $\hat \theta^n$ converge to $\hat \gamma$ and $\hat \theta$ respectively, and relabel this subsequence once more as $\{n\}$.
 Note that since $M^n$ solves $\clu(m^n, \eta^n)$, we have, for $t \in \RR_+$,
 $$M^n(t) = m^n - \int_0^t \eta^n(s) ds + \int_0^t M^n(s) ds.$$
 A straightforward calculation shows that
 $\|M^n(t) - M^n(s)\| \le 2 (t-s)$ for all $0 \le s \le t <\infty$,
 from which it follows that $\{M^n, \; n \in \NN \}$ is relatively compact in
 $C([0,\infty):\clp(\Delta^o))$. Assume without loss of generality (by selecting a further subsequence if needed) that $M^n \to M$ in 
 $C([0,\infty):\clp(\Delta^o))$ as $n\to \infty$.
 Note that we can write, for  $t \in \RR_+$ and $x \in \Delta^o$,
 $$M^n(t)(x) = m^n(x) - \int_0^t e^s \hat \gamma^n(\{x\} \times  ds) + \int_0^t M^n(s)(x) ds.$$
 Sending $n \to \infty$ in the previous display, we get
 \begin{equation}\label{eq:923}
 M(t)(x) = m(x) - \int_0^t e^s \hat \gamma(\{x\} \times  ds) + \int_0^t M(s)(x) ds.
 \end{equation}
 Furthermore, since $[\hat \gamma]_2(ds) = e^{-s} ds$, we can disintegrate $\hat \gamma$ as $\hat \gamma(\cdot \times ds) = \hat\eta(\cdot\mid s) e^{-s} ds$, where
 $s \mapsto \hat\eta(s) \doteq \hat \eta(\cdot \mid s)$ is a measurable map from $[0,\infty)$ to $\clp(\Delta^o)$. With this observation and \eqref{eq:923}, we have, for $t \in \RR_+$,
 $$
 M(t) = m - \int_0^t  \hat \eta(s) + \int_0^t M(s) ds
 $$
 which shows that 
 \begin{equation}\label{eq:940}
 \hat \eta \in \clu(m).
 \end{equation}
 
 Next, note that for each $t\in\RR_+$ and $x \in \Delta^o$
 $$\int_0^t e^s \hat \theta^n(\{x\}\times ds) = \int_0^t K(M^n(s))(x) ds.$$
 Sending $n \to \infty$ in the above display, we have, for each $ t \in \RR_+$ and $x \in \Delta^o$,
 $$\int_0^t e^s \hat \theta(\{x\}\times ds) = \int_0^t K(M(s))(x) ds.$$
 Also, since $[\hat \theta]_2(ds) = e^{-s} ds$, we can disintegrate $\hat \theta$ as $\hat \theta(\cdot \times ds) = \hat K(\cdot\mid s) e^{-s} ds$, where
 $s \mapsto \hat K(s) \doteq \hat K(\cdot \mid s)$ is a measurable map from $[0,\infty)$ to $\clp(\Delta^o)$. This says that for all $t \in \RR_+$ and $x \in \Delta^o$,
 $$\int_0^t K(M(s))(x) ds = \int_0^t \hat K(x \mid s) ds$$
 and so $K(M(s)) = \hat K(\cdot \mid s)$ for a.e. $s$.
 Finally, by the chain rule for relative entropies (see \cite[Corollary 2.7]{buddupbook})
 $$
 \int_0^{\infty} \exp(-s) R\left(\eta^n(\cdot \mid s) \| K(M^n(s)\right) ds = R(\hat \gamma^n \| \hat \theta^n),$$
 Using this, together with the lower semicontinuity of relative entropy
 and \eqref{eq:935} we now have
 \[
 R(\hat \gamma \| \hat \theta) \le
 \liminf_{n\to \infty} R(\hat \gamma^n \| \hat \theta^n) \le
 k.
 \]
 Using the chain rule again and the disintegrations of $\hat \gamma$ and $\hat \theta$ we now have
 \[
 \int_0^{\infty} \exp(-s) R\left(\hat\eta(\cdot \mid s) \| K(M(s)\right) ds  = R(\hat \gamma \| \hat \theta) \le k.
 \]
 Together this, \eqref{eq:940}, and the definition of $I$ show that
 $I(m) \le k$. The result follows. \hfill \qed
 $\,$\\
 
 \noindent {\bf Acknowledgements.} Research of AB supported in part by  the NSF (DMS-1814894, DMS-1853968 and DMS-2134107).
 \bibliographystyle{cas-model2-names}
\bibliography{cas-refs}

\end{document}